# Dual-Thrust Switching Analytical Guidance Algorithm for Powered Landing with Attitude Smoothness Optimization

Wen-Bo Li[1]

*School of Aerospace Engineering, Tsinghua University, Beijing 100084, China*

Dai Shen[2]

*School of Astronautics, Beihang University, Beijing, 102206, China*

Sheng-Ping Gong[3*]

*School of Astronautics, Beihang University, Beijing, 102206, China*

*State Key Laboratory of High-Efficiency Reusable Aerospace Transportation Technology, Beijing, 102206, China*

**Traditional numerical guidance methods for powered landing of reusable rockets are typically constrained by high computational complexity and inadequate real-time performance. Moreover, insufficient consideration of attitude smoothness often induces severe fluctuations in control commands; meanwhile, most existing approaches are tailored for single-thrust scenarios, failing to accommodate the guidance requirements of multi-engine thrust switching. To mitigate these limitations, this paper proposes an analytical guidance method optimized for attitude smoothness, which supports dual-thrust-mode switching. First, a corresponding optimal control problem is formulated, and it is theoretically proven that the optimal attitude command takes a concise piecewise cubic function form. This transforms complex trajectory optimization into a parametric analytical optimization problem, yielding a substantial improvement in computational efficiency. Further, a three-phase guidance framework is designed to enable adaptive determination of the guidance activation point and thrust switching point; when integrated with an aerodynamic correction strategy, this framework enhances the method's adaptability in complex flight environments—particularly under high lift-to-drag ratio conditions. Simulation results demonstrate that the attitude command profile generated by the proposed method aligns closely with the theoretical optimal solution, with an ultra-short computation time, confirming its strong potential for online real-time implementation. Even**

1 Ph.D. Candidate, School of Aerospace Engineering, Tsinghua University. Email: liwb21@mails.tsinghua.edu.cn

2 Master's Candidate, School of Astronautics, Beihang University.

3* Professor, School of Astronautics, Beihang University. Email: gongsp@buaa.edu.cn (corresponding author)

**under stringent conditions (e.g., limited thrust adjustment range, high lift-to-drag ratios, and parameter deviations), the method consistently achieves high-precision landing, showcasing promising prospects for engineering applications.**

## I. Introduction

Reusable launch vehicles (RLVs) represent one of the key technologies for achieving low-cost and high-frequency access to space, with one of their core procedures being the realization of high-precision and high-reliability powered landing during the return phase. The design of powered landing guidance laws is a typical multi-objective and strongly constrained trajectory optimization problem [1]. It must not only satisfy the stringent constraints on terminal position, velocity, and attitude but also balance multiple requirements including the stability of the control system, the thrust adjustment capability of engines, and the real-time performance of computations.

Early methods such as gravity-turn guidance [2, 3] and explicit polynomial guidance [4, 5] exhibit strong engineering practicability; however, they face significant limitations in incorporating aerodynamic forces and process constraints, and the algorithms themselves lack optimality—restricting their application in recoverable rockets. Over the past two decades, research focus has gradually shifted toward numerical solution methods for optimal control, which are primarily categorized into indirect methods and direct methods. Indirect methods apply the Pontryagin's Maximum Principle to transform the original problem into a two-point boundary value problem (TPBVP) [4, 6-9]. While these methods guarantee the optimality of solutions, they are highly sensitive to the initial values of costate variables and have difficulty dealing with path constraints, thus limiting their practical utility. Among direct methods, convex optimization techniques have garnered considerable attention due to their deterministic advantage of convergence in polynomial time [10]. The "lossless convexification" technique proposed by Acikmese and Ploen [11] successfully addresses the lower thrust bound constraint and provides a globally optimal solution for the fuel-optimal landing problem. Since then, this technique has been extended to handle complex scenarios such as thrust pointing constraints [12-14], affine and quadratic state constraints [15-17], and nonlinear mixed-integer dynamic models [18, 19]. Sequential convex optimization methods, including SCvx [20, 21] and GuSTO [22], further solve strongly nonlinear and nonconvex problems (e.g., aerodynamic forces [23] and six-degree-of-freedom (6-DOF) models [24]) by iteratively solving local convex approximations. In addition, numerous scholars have investigated convex optimization methods that incorporate stochastic terms into dynamic models [25], and the introduction of chance constraints has significantly expanded the engineering application scope of such algorithms [26-28].

Traditional powered landing guidance methods largely rely on numerical optimization techniques such as sequential convex optimization and pseudospectral methods. Although these methods can handle complex dynamic models and constraints, they suffer from high computational burden, convergence dependence on initial guesses, and difficulty in ensuring real-time performance—limiting their direct implementation on on-board computers. Moreover, existing methods often prioritize fuel optimality while neglecting the smoothness of attitude angle variations, which may lead to severe fluctuations in attitude and thrust during actual flight and reduce system reliability.

Meanwhile, most existing research focuses on the guidance problem under the single-thrust mode, whereas engineering practice has demonstrated that the dual-thrust mode offers greater application value. The dual-thrust mode refers to a strategy in the powered landing process where multiple engines are ignited simultaneously for rapid deceleration at the initial stage; subsequently, according to preset switching conditions, some engines are shut down, and the system switches to low-thrust mode for precise maneuvering, thereby achieving a smooth landing. Currently, this strategy is adopted by vehicles such as the first stage of Starship and New Glenn. Compared with the single-thrust mode, the dual-thrust mode presents two main advantages:

(1) It can balance the relationship between fuel consumption and landing thrust-to-weight ratio: rapid deceleration at the initial stage helps reduce gravity losses and ensure fuel economy, while smooth landing at the terminal stage improves terminal accuracy and increases the feasible region [29].

(2) It shows stronger adaptability to aerodynamic parameter uncertainties: these uncertainties are typically significant and greatly affect guidance, especially during the initial phase of powered landing. In complex aerodynamic environments, the single-thrust mode relies on engines' deep thrust modulation capability to achieve precise landing. By contrast, the dual-thrust mode partially compensates for aerodynamic disturbances through adaptive adjustment of thrust switching time, thus reducing requirements for engine thrust modulation.

To overcome the challenges of current research and fully leverage the advantages of the dual-thrust mode, this paper proposes an analytical guidance algorithm with attitude smoothness optimization, named SAM-DS (Smooth Attitude Maneuvering with Dual-thrust Switching). Through rigorous mathematical derivation, this method constructs a closed-loop analytical solution, transforming the complex trajectory optimization problem into a series of analytically solvable parameter optimization problems. This fundamentally eliminates the need for iterative optimization, significantly improving computational efficiency—the calculation time is reduced by 2 to 3 orders of magnitude compared with traditional numerical methods—thus facilitating online real-time execution. Open-loop simulation results indicate that the output attitude angle profile is highly consistent with the theoretical optimal solution, verifying the optimality of the algorithm in terms of attitude smoothness. The main theoretical innovations of this paper can be summarized as follows:

(1) An analytical optimal control problem with the objective of attitude angle smoothness is established. It is proven that the optimal solution is a piecewise cubic function of time, satisfying continuity, differentiability, and proportionality conditions, which lays a theoretical foundation for the optimality of the algorithm.

(2) A three-phase guidance framework consisting of the "pre-activation phase", "high-thrust phase", and "low-thrust phase" is proposed, along with an adaptive transition criterion for the activation state and thrust switching state. This significantly reduces the algorithm's requirement for the engine's deep thrust modulation capability.

(3) An aerodynamic acceleration correction strategy based on piecewise polynomial interpolation is proposed, which greatly enhances the robustness and adaptability of the algorithm in complex environments, especially in scenarios with high lift-to-drag ratios and significant initial aerodynamic forces.

In conclusion, the proposed algorithm features a clear process and low parameter dependence, demonstrating excellent interpretability and portability. It is suitable for practical engineering systems with strict requirements on computational resources and real-time performance.

## II. Dynamics Modeling and Simplification

This study focuses on the powered landing phase of a recoverable rocket. A flat Earth model is adopted, the non-inertial force induced by Earth’s rotation is neglected, and the aerodynamic force is calculated using an ellipsoid model. First, two coordinate systems are defined: the landing coordinate system $O_l - x_l y_l z_l$ and the body coordinate system $O_b - x_b y_b z_b$ . In the landing coordinate system, the origin $O_l$ is located at the landing site; the $x_l$ -axis points along the flight path and is parallel to the ground plane; the $y_l$ -axis points to the zenith;

and the $z_l$ -axis is determined by the right-hand rule. In the body coordinate system, the origin $O_b$ is located at the rocket's center of mass; the $x_b$ -axis is along the rocket's longitudinal axis pointing to the nose; the $y_b$ -axis lies in the rocket's longitudinal plane; and the $z_b$ -axis is determined by the right-hand rule.

The projection of the three-degree-of-freedom (3-DOF) dynamic equations in the landing coordinate system can be expressed as follows:

$$\dot{\boldsymbol{r}} = \boldsymbol{v},\ \dot{\boldsymbol{v}} = \boldsymbol{a}_T + \boldsymbol{a}_R + \boldsymbol{g},\ \dot{m} = -\frac{T}{I_{\text{sp}}}, \boldsymbol{a}_T = \frac{T}{m}\boldsymbol{L}_{bl}^T\begin{bmatrix}1\\0\\0\end{bmatrix},\ \boldsymbol{a}_R = -\frac{\rho S\|\boldsymbol{v}\|}{2m}\boldsymbol{L}_{bl}^T\boldsymbol{C}_R\boldsymbol{L}_{bl}\boldsymbol{v} \tag{1}$$

where $\boldsymbol{r} = [r_x, r_y, r_z]^T$ denotes the position projection in the landing coordinate system, $\boldsymbol{v} = [v_x, v_y, v_z]^T$ denotes the velocity projection in the landing coordinate system, $m$ is the rocket mass, $\boldsymbol{a}_T$ is the thrust acceleration projection, $\boldsymbol{a}_R$ is the aerodynamic acceleration projection, $\boldsymbol{g} = [0, -g, 0]^T$ is the gravitational acceleration projection, $\boldsymbol{L}_{bl}$ is the transformation matrix from the landing coordinate system to the body coordinate system, $I_{\text{sp}}$ is the specific impulse of the engine, $T$ is the thrust magnitude, $\rho$ is the constant atmospheric density, $S$ is the reference area of the rocket, $\|\bullet\|$ represents the 2-norm of a vector, and $\boldsymbol{C}_R$ is the ellipsoid aerodynamic coefficient matrix. $\boldsymbol{C}_R$ and $\boldsymbol{L}_{bl}$ can be expanded as follows:

$$\boldsymbol{C}_R = \begin{bmatrix}C_x & 0 & 0\\0 & C_y & 0\\0 & 0 & C_z\end{bmatrix},\ \boldsymbol{L}_{bl} = \begin{bmatrix}\cos\psi & 0 & -\sin\psi\\0 & 1 & 0\\\sin\psi & 0 & \cos\psi\end{bmatrix}\begin{bmatrix}\cos\theta & \sin\theta & 0\\-\sin\theta & \cos\theta & 0\\0 & 0 & 1\end{bmatrix} \tag{2}$$

where $\psi$ denotes the yaw angle and $\theta$ denotes the pitch angle. By treating the thrust as a piecewise constant value for guidance law design, deep thrust modulation is avoided, and the differential equation describing the rocket mass variation can be directly eliminated. The dynamic model can be further rewritten as follows:

$$\begin{gathered}
\dot{\boldsymbol{r}} = \boldsymbol{v},\ \dot{\boldsymbol{v}} = \boldsymbol{a}_T(T,\theta,\psi,\tau) + \boldsymbol{a}_R(\theta,\psi,\boldsymbol{v},\tau) + \boldsymbol{g},\ \dot{T} = 0\\
\boldsymbol{a}_T(T,\theta,\psi,\tau) = \begin{cases}\dfrac{nT}{m_0 - \dfrac{nT}{I_{\text{sp}}}\tau}\boldsymbol{L}_{bl}^T\begin{bmatrix}1\\0\\0\end{bmatrix},\ \tau \le \mu t_{\text{go}}^{(2)}\\[2ex]\dfrac{T}{m_0 - \dfrac{nT}{I_{\text{sp}}}\mu t_{\text{go}}^{(2)} - \dfrac{T}{I_{\text{sp}}}\left(\tau - \mu t_{\text{go}}^{(2)}\right)}\boldsymbol{L}_{bl}^T\begin{bmatrix}1\\0\\0\end{bmatrix},\ \mu t_{\text{go}}^{(2)} < \tau \le (1+\mu)\, t_{\text{go}}^{(2)}\end{cases}\\
\boldsymbol{a}_R(\theta,\psi,\boldsymbol{v},\tau) = \begin{cases}-\dfrac{\rho S\|\boldsymbol{v}\|}{2\left(m_0 - \dfrac{nT}{I_{\text{sp}}}\tau\right)}\boldsymbol{L}_{bl}^T\boldsymbol{C}_R\boldsymbol{L}_{bl}\begin{bmatrix}v_x\\v_y\\v_z\end{bmatrix},\ \tau \le \mu t_{\text{go}}^{(2)}\\[2ex]-\dfrac{\rho S\|\boldsymbol{v}\|}{2\left(m_0 - \dfrac{nT}{I_{\text{sp}}}\mu t_{\text{go}}^{(2)} - \dfrac{T}{I_{\text{sp}}}\left(\tau - \mu t_{\text{go}}^{(2)}\right)\right)}\boldsymbol{L}_{bl}^T\boldsymbol{C}_R\boldsymbol{L}_{bl}\begin{bmatrix}v_x\\v_y\\v_z\end{bmatrix},\ \mu t_{\text{go}}^{(2)} < \tau \le (1+\mu)\, t_{\text{go}}^{(2)}\end{cases}
\end{gathered} \tag{3}$$

where $m_0$ is the initial rocket mass and $\tau$ is the flight time. In this expression, $n$ denotes the thrust ratio coefficient between the high-thrust and low-thrust phases, which is a known constant greater than 1. $t_{\text{go}}^{(1)}$ and $t_{\text{go}}^{(2)}$ represent the flight times of the high-thrust and low-thrust phases, respectively, and the total flight time is given by $t_{\text{go}} = t_{\text{go}}^{(1)} + t_{\text{go}}^{(2)}$. The parameter $\mu = t_{\text{go}}^{(1)} / t_{\text{go}}^{(2)}$ is referred to as the time allocation coefficient. The dynamic equation $\dot{T} = 0$ is supplemented herein. It is worth noting that $\dot{T} = 0$ indicates that the thrust magnitude is constant over time, although the specific value of the thrust remains undetermined. Next, the dynamic model is simplified. Considering that the rocket essentially moves within a longitudinal plane (the x-y plane), the yaw angle can be treated as a small quantity. Substituting Eq. (2) into Eq. (3) and utilizing a Taylor expansion, $\boldsymbol{a}_T\left(T,\theta,\psi,\tau\right)$ can be approximated as follows:

$$\begin{aligned}
&\boldsymbol{a}_T\left(T,\theta,\psi,\tau\right) \approx \tilde{\boldsymbol{a}}_T\left(T,\Delta\theta,\Delta\psi,\tau\right) \\
&= \begin{cases}
\dfrac{nT}{m_0}\begin{bmatrix} \cos\bar{\theta} - \sin\bar{\theta}\Delta\theta - \dfrac{1}{2}\cos\bar{\theta}\Delta\theta^2 \\ \sin\bar{\theta} + \cos\bar{\theta}\Delta\theta - \dfrac{1}{2}\sin\bar{\theta}\Delta\theta^2 \\ -\left(\bar{\psi} + \Delta\psi\right) \end{bmatrix}\left(1 + \dfrac{nT}{I_{\text{sp}}m_0}\tau\right), & \tau \le \mu t_{\text{go}}^{(2)} \\
\dfrac{T}{m_0}\begin{bmatrix} \cos\bar{\theta} - \sin\bar{\theta}\Delta\theta - \dfrac{1}{2}\cos\bar{\theta}\Delta\theta^2 \\ \sin\bar{\theta} + \cos\bar{\theta}\Delta\theta - \dfrac{1}{2}\sin\bar{\theta}\Delta\theta^2 \\ -\left(\bar{\psi} + \Delta\psi\right) \end{bmatrix}\left(1 + \dfrac{(n-1)T}{I_{\text{sp}}m_0}\mu t_{\text{go}}^{(2)} + \dfrac{T}{I_{\text{sp}}m_0}\tau\right), & \mu t_{\text{go}}^{(2)} < \tau \le \left(1+\mu\right)t_{\text{go}}^{(2)}
\end{cases}
\end{aligned} \tag{4}$$

where $\bar{\theta}$ and $\bar{\psi}$ are undetermined constant principal components, and $\Delta\theta$ and $\Delta\psi$ are undetermined functions of time. If the influence of mass variation on the aerodynamic acceleration is neglected, substituting Eq. (2) into Eq. (3) also simplifies $\boldsymbol{a}_R\left(\theta,\psi,\boldsymbol{v},\tau\right)$ to:

$$\boldsymbol{a}_R\left(\theta,\psi,\boldsymbol{v},\tau\right) \approx \tilde{\boldsymbol{a}}_R\left(\theta,\boldsymbol{v}\right) = -k\left\|\boldsymbol{v}\right\|\begin{bmatrix} f_x(\theta,v_x,v_y) \\ f_y(\theta,v_x,v_y) \\ f_z(v_z) \end{bmatrix} \tag{5}$$

where

$$\begin{aligned}
&f_x(\theta,v_x,v_y) = \left(\cos^2\theta + \frac{C_y}{C_x}\sin^2\theta\right)v_x + \left(1 - \frac{C_y}{C_x}\right)\cos\theta\sin\theta v_y \\
&f_y(\theta,v_x,v_y) = \left(\sin^2\theta + \frac{C_y}{C_x}\cos^2\theta\right)v_y + \left(1 - \frac{C_y}{C_x}\right)\cos\theta\sin\theta v_x \\
&f_z(v_z) = \frac{C_z}{C_x}v_z \\
&k = \frac{\rho S C_x}{2m_0}
\end{aligned} \tag{6}$$

To facilitate analytical integration, the aerodynamic acceleration is decomposed into a principal component and an incremental component, both of which are approximated as functions of time:

$$\tilde{\boldsymbol{a}}_R(\theta,\boldsymbol{v}) \approx \bar{\boldsymbol{a}}_R(\tau)+\Delta\boldsymbol{a}_R(\tau)$$
$$\bar{\boldsymbol{a}}_R(\tau)=\begin{cases}-k\left\|\boldsymbol{v}_0\right\|\left[1-\dfrac{n\tau}{(n\mu+1)t_{\mathrm{go}}^{(2)}}\right]^2\boldsymbol{v}_0,\ \tau\le\mu t_{\mathrm{go}}^{(2)}\\ -k\left\|\boldsymbol{v}_0\right\|\left[\dfrac{(\mu+1)t_{\mathrm{go}}^{(2)}-\tau}{(n\mu+1)t_{\mathrm{go}}^{(2)}}\right]^2\boldsymbol{v}_0,\ \mu t_{\mathrm{go}}^{(2)}<\tau\le(1+\mu)t_{\mathrm{go}}^{(2)}\end{cases} \tag{7}$$
$$\Delta\boldsymbol{a}_R(\tau)=\begin{cases}\hat{\boldsymbol{a}}_R^{(1)}\tau^3+\hat{\boldsymbol{b}}_R^{(1)}\tau^2+\hat{\boldsymbol{c}}_R^{(1)}\tau+\hat{\boldsymbol{d}}_R^{(1)},\ \tau\le\mu t_{\mathrm{go}}^{(2)}\\ \hat{\boldsymbol{a}}_R^{(2)}\left(\tau-t_{\mathrm{go}}^{(1)}\right)^3+\hat{\boldsymbol{b}}_R^{(2)}\left(\tau-t_{\mathrm{go}}^{(1)}\right)^2+\hat{\boldsymbol{c}}_R^{(2)}\left(\tau-t_{\mathrm{go}}^{(1)}\right)+\hat{\boldsymbol{d}}_R^{(2)},\ \mu t_{\mathrm{go}}^{(2)}<\tau\le(1+\mu)t_{\mathrm{go}}^{(2)}\end{cases}$$

where $\boldsymbol{v}_0=[v_{x0},v_{y0},v_{z0}]^T$ is the initial velocity projection of the rocket, and $\bar{\boldsymbol{a}}_R(\tau)$ characterizes the aerodynamic acceleration under the assumption of piecewise linear velocity decrease. $\Delta\boldsymbol{a}_R(\tau)$ represents the nonlinear component, and $\hat{\boldsymbol{a}}_R^{(i)}\sim\hat{\boldsymbol{d}}_R^{(i)}$ are interpolation coefficients. It is worth noting that when the rocket's lift-to-drag ratio and initial velocity are small, $\Delta\boldsymbol{a}_R(\tau)$ can be treated as a small quantity; otherwise, it cannot be neglected and must be estimated through a specific method. This work proposes a piecewise polynomial interpolation method for this estimation, which exhibits excellent real-time performance (see Section V.F for details). Consequently, based on Eqs. (4) and (7), the dynamic equation (3) can be further simplified as follows:

$$\dot{\boldsymbol{r}}=\boldsymbol{v},\ \dot{\boldsymbol{v}}=\tilde{\boldsymbol{a}}_T(T,\Delta\theta,\Delta\psi,\tau)+\bar{\boldsymbol{a}}_R(\tau)+\Delta\boldsymbol{a}_R(\tau)+\boldsymbol{g},\ \dot{T}=0 \tag{8}$$

Subsequent derivations will be based on the dynamic model given by Eq.(8).

## III. Overall Scheme Design

Based on the different treatments of terminal constraints, this paper first defines two types of problems: the "fixed-point landing problem (FLP)" and the "quasi-fixed-point landing problem (QLP)", laying the foundation for the subsequent construction of the analytical guidance law.

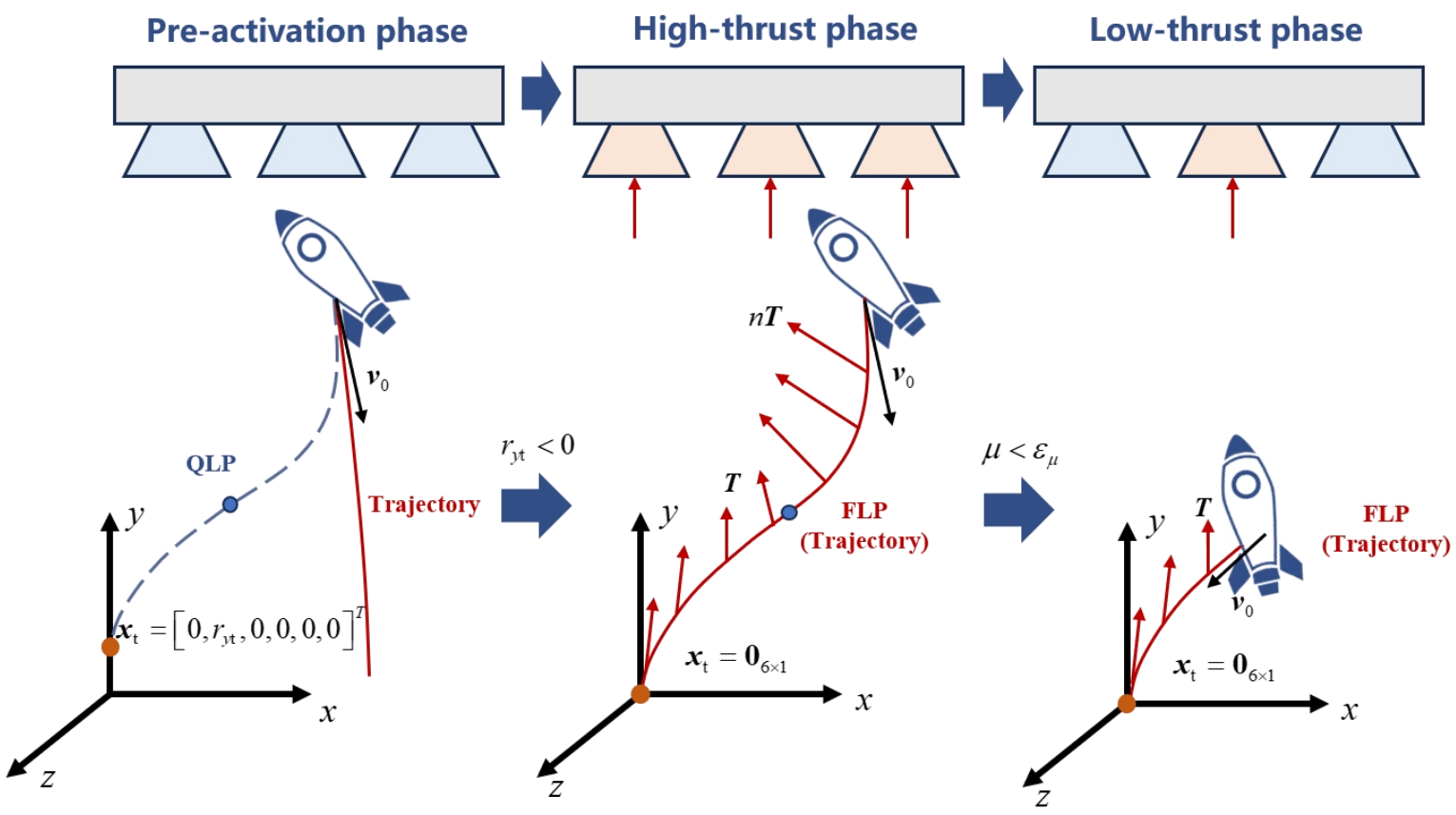

**Fig. 1 Schematic diagram of the landing process**

Among them, the target state constraint for the FLP is $\boldsymbol{x}_{\mathrm{t}}=\mathbf{0}_{6\times1}$ . The solution to this problem enables the rocket to achieve a soft landing at the target point. The QLP assumes that the rocket's thrust value $T$ is a given constant function (nominal thrust), does not constrain the terminal height, and only requires the rocket to softly land "above" the specified point. The target state constraint is $\boldsymbol{x}_{\mathrm{t}}=[0,r_{y\mathrm{t}},0,0,0,0]^{T}$ , where $r_{y\mathrm{t}}$ is a free variable. After solving this problem, the sign of $r_{y\mathrm{t}}$ can be used to determine whether the rocket should activate guidance. A schematic diagram is shown in Fig. 1.

In view of the above analysis, the powered landing can be divided into three phases, referred to as the "pre-activation phase," the "high-thrust phase," and the "low-thrust phase." In the pre-activation phase, the predicted terminal height $r_{y\mathrm{t}}$ is output by solving the QLP in each guidance cycle. When $r_{y\mathrm{t}}\le0$ , the rocket enters the high-thrust phase. In each guidance cycle, the FLP is solved to output the guidance commands pitch angle $\theta$ and yaw angle $\psi$ , and fine-tune the time allocation coefficient $\mu=t_{\mathrm{go}}^{(1)}/t_{\mathrm{go}}^{(2)}$ until $\mu\le\varepsilon_{\mu}$ (where $\varepsilon_{\mu}$ is a small threshold). Subsequently, the rocket enters the low-thrust phase. In each guidance cycle, the FLP is solved analytically to output the guidance commands pitch angle $\theta$ and yaw angle $\psi$ , and fine-tune the thrust magnitude $T$, ultimately achieving landing. The following is the flowchart of the overall scheme.

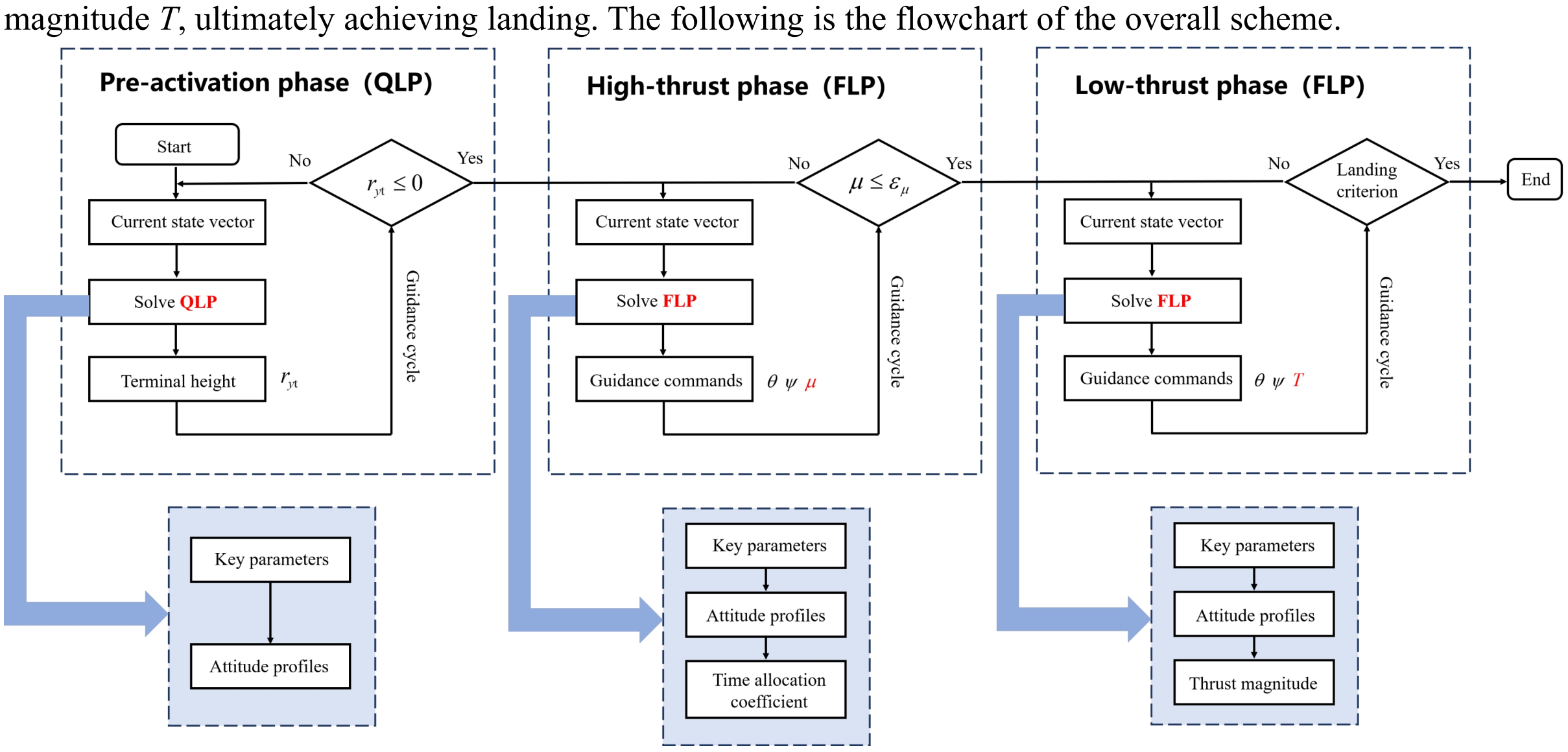

**Fig. 2 Schematic diagram of the guidance algorithm flow**

As can be seen from Fig. 2, deriving the analytical solution for the FLP and QLP presents the primary challenge and constitutes the core the algorithm. One of the main contributions of this paper is to propose an analytical method, where the output guidance commands achieve the optimal smoothness of attitude angle. Subsequently, the FLP is mathematically formulated as Problem 1.

**Problem 1:**

$$\underset{u1,u2,u3}{\text{Minimize}}\ J=\int_{\tau=0}^{(1+\mu)t_{\text{go}}^{(2)}}\frac{1}{2}\left(u_1^2+u_2^2\right)\text{d}\,\tau$$

$$\dot{\boldsymbol{r}}=\boldsymbol{v},\ \dot{\boldsymbol{v}}=\tilde{\boldsymbol{a}}_T\left(T,\Delta\theta,\Delta\psi,\tau\right)+\bar{\boldsymbol{a}}_R\left(\tau\right)+\Delta\boldsymbol{a}_R\left(\tau\right)+\boldsymbol{g},\ \dot{T}=0$$

$$\Delta\dot{\theta}=u_1,\ \Delta\dot{\psi}=u_2,\ \dot{T}=u_3,\ u_3=0$$

$$\left[\boldsymbol{r}^T,\boldsymbol{v}^T\right]^T\Big|_{\tau=0}=\left[\boldsymbol{r}_0^T,\boldsymbol{v}_0^T\right]^T,\left[\boldsymbol{r}^T,\boldsymbol{v}^T\right]^T\Big|_{\tau=(1+\mu)t_{\text{go}}^{(2)}}=\boldsymbol{0}_{6\times1}$$

$$\Delta\theta\Big|_{\tau=(1+\mu)t_{\text{go}}^{(2)}}=\frac{\pi}{2}-\bar{\theta},\Delta\psi\Big|_{\tau=(1+\mu)t_{\text{go}}^{(2)}}=-\bar{\psi}$$

To resolve this problem, Theorem 1 is presented.

**Theorem 1**: If $\frac{(n-1)T}{I_{\text{sp}}m_0}\mu t_{\text{go}}^{(2)}+\frac{T}{I_{\text{sp}}m_0}\tau\ll1$, $\left|\frac{1}{2}\cos\bar{\theta}\Delta\theta^2\right|\ll1$ and $\left|\frac{1}{2}\sin\bar{\theta}\Delta\theta^2\right|\ll1$, the optimal solutions $\Delta\theta^*$ and $\Delta\psi^*$ to Problem 1 are piecewise cubic functions of time:

$$\begin{aligned}
\Delta\theta^*&=\begin{cases}a^{(1)}\left(\tau-t_{\text{go}}\right)^3+b^{(1)}\left(\tau-t_{\text{go}}\right)^2+c^{(1)}\left(\tau-t_{\text{go}}\right)+d^{(1)},\tau\le\mu t_{\text{go}}^{(2)}\\ a^{(2)}\left(\tau-t_{\text{go}}\right)^3+b^{(2)}\left(\tau-t_{\text{go}}\right)^2+c^{(2)}\left(\tau-t_{\text{go}}\right)+d^{(2)},\mu t_{\text{go}}^{(2)}<\tau\le(\mu+1)t_{\text{go}}^{(2)}\end{cases}\\
\Delta\psi^*&=\begin{cases}a'^{(1)}\left(\tau-t_{\text{go}}\right)^3+b'^{(1)}\left(\tau-t_{\text{go}}\right)^2+c'^{(1)}\left(\tau-t_{\text{go}}\right)+d'^{(1)},\tau\le\mu t_{\text{go}}^{(2)}\\ a'^{(2)}\left(\tau-t_{\text{go}}\right)^3+b'^{(2)}\left(\tau-t_{\text{go}}\right)^2+c'^{(2)}\left(\tau-t_{\text{go}}\right)+d'^{(2)},\mu t_{\text{go}}^{(2)}<\tau\le(\mu+1)t_{\text{go}}^{(2)}\end{cases}
\end{aligned}\tag{9}$$

This piecewise cubic function satisfies the following junction conditions:

(1) Continuity Condition: $\Delta\theta^*\Big|_{\tau=t_{\text{go}}^{(1)-}}=\Delta\theta^*\Big|_{\tau=t_{\text{go}}^{(1)+}},\Delta\psi^*\Big|_{\tau=t_{\text{go}}^{(1)-}}=\Delta\psi^*\Big|_{\tau=t_{\text{go}}^{(1)+}}$

(2) Differentiability Condition: $\frac{\text{d}\Delta\varphi^*}{\text{d}\tau}\Big|_{\tau=t_{\text{go}}^{(1)-}}=\frac{\text{d}\Delta\varphi^*}{\text{d}\tau}\Big|_{\tau=t_{\text{go}}^{(1)+}},\frac{\text{d}\Delta\psi^*}{\text{d}\tau}\Big|_{\tau=t_{\text{go}}^{(1)-}}=\frac{\text{d}\Delta\psi^*}{\text{d}\tau}\Big|_{\tau=t_{\text{go}}^{(1)+}}$

(3) Proportionality Condition: $a^{(1)}=na^{(2)},b^{(1)}=nb^{(2)},a'^{(1)}=na'^{(2)},b'^{(1)}=nb'^{(2)}$

The detailed proof is provided in Appendix A.

It is worth noting that Theorem 1 is also applicable to the QLP. Theorem 1 indicates that the optimal solution to Problem 1 exhibits favorable properties, simplifying the complex optimal control problem into a parameter optimization problem involving 8 undetermined coefficients of piecewise cubic functions. This effectively avoids the need for solving complex numerical optimization problems.

In summary, the specific methodology for solving both the FLP and the QLP is outlined as follows:

(1) Calculation of Key Parameters:

First, based on Eq. (8), the principal components of the pitch angle $\bar{\theta}$ and yaw angle $\bar{\psi}$, as well as the remaining flight times $t_{\text{go}}^{(1)}$ and $t_{\text{go}}^{(2)}$, are derived to satisfy the terminal velocity constraints. The detailed implementation is presented in Section V.A.

(2) Calculation of Attitude Angle Profiles:

Subsequently, $\Delta\theta$ and $\Delta\psi$ are designed as piecewise cubic polynomial functions with undetermined coefficients. The coefficients of these polynomials (16 parameters in total for both pitch and yaw channels) are computed analytically using constraints such as horizontal terminal position, terminal vertical attitude, the objective function, and junction conditions. This derivation yields the current pitch angle increment $\Delta\theta$ and yaw angle increment $\Delta\psi$. The detailed implementation is presented in Section V.B and Section V.C.

(3) Calculation of Time Allocation Coefficient:

For the high-thrust phase, to ensure the terminal height is zero, the time allocation ratio $\mu$ is fine-tuned according to the terminal height constraint. The detailed implementation is presented in Section V.D.

(4) Calculation of Thrust Magnitude:

For the low-thrust phase, to ensure the terminal height is zero, the thrust magnitude $T$ is fine-tuned according to the terminal height constraint. The detailed implementation is presented in Section V.E.

## IV. Analytical Representation of Constraints

Prior to deriving the guidance algorithm, it is necessary to represent the terminal constraints using analytical methods. A series of integral functionals are defined next to facilitate the representation of these constraints. When $\Delta\theta$ and $\Delta\psi$ are piecewise cubic functions, the following integral functionals can all be expressed in explicit form. The specific expressions are given by Eqs. (B5), (B7), (B11), (B13), (B14), and (B15) in Appendix B.

$$\begin{aligned}
\Theta(x) &= n\Theta^{(1)}(x)+\Theta^{(2)}(x)=n\int_{\tau=0}^{\mu t_{\mathrm{go}}^{(2)}} x\left(1+nh\tau\right)\,\mathrm{d}\tau+\int_{\tau=\mu t_{\mathrm{go}}^{(2)}}^{(\mu+1)t_{\mathrm{go}}^{(2)}} x\left[1+(n-1)h\mu t_{\mathrm{go}}^{(2)}+h\tau\right]\mathrm{d}\tau \\
\tilde{\Theta}(x) &= nt_{\mathrm{go}}^{(2)}\Theta^{(1)}(x)+n\tilde{\Theta}^{(1)}(x)+\tilde{\Theta}^{(2)}(x) \\
&= nt_{\mathrm{go}}^{(2)}\int_{\tau=0}^{\mu t_{\mathrm{go}}^{(2)}} x\left(1+nh\tau\right)\,\mathrm{d}\tau+n\int_{t=0}^{\mu t_{\mathrm{go}}^{(2)}}\int_{\tau=0}^{t} x\left(1+nh\tau\right)\,\mathrm{d}\tau\mathrm{d}t+\int_{t=\mu t_{\mathrm{go}}^{(2)}}^{(\mu+1)t_{\mathrm{go}}^{(2)}}\int_{\tau=\mu t_{\mathrm{go}}^{(2)}}^{t} x\left[1+(n-1)h\mu t_{\mathrm{go}}^{(2)}+h\tau\right]\mathrm{d}\tau\mathrm{d}t
\end{aligned} \tag{10}$$

where $h=T/(I_{\mathrm{sp}}m_0)$, the superscript denotes the phase number, and the tilde symbol represents double integration. Directly integrating Eq. (8) and substituting Eq. (10), the velocity components at the thrust switching point can be obtained as:

$$\boldsymbol{v}_{\mathrm{m}}=\boldsymbol{v}_0+\frac{nT}{m_0}\begin{bmatrix}\cos\bar{\theta}\left[\chi_{\mathrm{m}}-\frac{1}{2}\Theta^{(1)}\left(\Delta\theta^2\right)\right]-\Theta^{(1)}\left(\Delta\theta\right)\sin\bar{\theta} \\ \sin\bar{\theta}\left[\chi_{\mathrm{m}}-\frac{1}{2}\Theta^{(1)}\left(\Delta\theta^2\right)\right]+\Theta^{(1)}\left(\Delta\theta\right)\cos\bar{\theta} \\ -\bar{\psi}\chi_{\mathrm{m}}-\Theta^{(1)}\left(\Delta\psi\right)\end{bmatrix}+\mu t_{\mathrm{go}}^{(2)}\boldsymbol{g}-k\left\|\boldsymbol{v}_0\right\|t_{\mathrm{go}}^{(2)}\Pi_{\mathrm{m}}\boldsymbol{v}_0+\boldsymbol{\xi}_{\mathrm{m}} \tag{11}$$

where

$$
\begin{aligned}
\Pi_{\mathrm{m}} &= \frac{n^2\mu^3 + 3n\mu^2 + 3\mu}{3\left(n\mu+1\right)^2} \\
\chi_{\mathrm{m}} &= \mu t_{\mathrm{go}}^{(2)} + \frac{nT}{2I_{\mathrm{sp}}m_0}\mu^2 t_{\mathrm{go}}^{(2)2} \\
\boldsymbol{\xi}_{\mathrm{m}} &= \frac{1}{4}\hat{\boldsymbol{a}}_R^{(1)}\left(\mu t_{\mathrm{go}}^{(2)}\right)^4 + \frac{1}{3}\hat{\boldsymbol{b}}_R^{(1)}\left(\mu t_{\mathrm{go}}^{(2)}\right)^3 + \frac{1}{2}\hat{\boldsymbol{c}}_R^{(1)}\left(\mu t_{\mathrm{go}}^{(2)}\right)^2 + \hat{\boldsymbol{d}}_R^{(1)}\mu t_{\mathrm{go}}^{(2)}
\end{aligned}
\tag{12}
$$

Similarly, the expressions for the terminal velocity component constraints are:

$$
\mathbf{0} = \boldsymbol{v}_0 + \frac{T}{m_0}\begin{bmatrix} \cos\bar{\theta}\left[\chi - \frac{1}{2}\Theta\left(\Delta\theta^2\right)\right] - \Theta\left(\Delta\theta\right)\sin\bar{\theta} \\ \sin\bar{\theta}\left[\chi - \frac{1}{2}\Theta\left(\Delta\theta^2\right)\right] + \Theta\left(\Delta\theta\right)\cos\bar{\theta} \\ -\bar{\psi}\chi - \Theta\left(\Delta\psi\right) \end{bmatrix} + \left(\mu+1\right)t_{\mathrm{go}}^{(2)}\boldsymbol{g} - k\left\|\boldsymbol{v}_0\right\|t_{\mathrm{go}}^{(2)}\Pi\boldsymbol{v}_0 + \boldsymbol{\xi}
\tag{13}
$$

where

$$
\begin{aligned}
\Pi &= \frac{n^2\mu^3 + 3n\mu^2 + 3\mu + 1}{3\left(n\mu+1\right)^2} \\
\chi &= \left(n\mu+1\right)t_{\mathrm{go}}^{(2)} + \frac{\left(n\mu+1\right)^2 T}{2I_{\mathrm{sp}}m_0}t_{\mathrm{go}}^{(2)2} \\
\boldsymbol{\xi} &= \frac{1}{4}\left(\hat{\boldsymbol{a}}_R^{(2)} + \hat{\boldsymbol{a}}_R^{(1)}\mu^4\right)t_{\mathrm{go}}^{(2)4} + \frac{1}{3}\left(\hat{\boldsymbol{b}}_R^{(2)} + \hat{\boldsymbol{b}}_R^{(1)}\mu^3\right)t_{\mathrm{go}}^{(2)3} + \frac{1}{2}\left(\hat{\boldsymbol{c}}_R^{(2)} + \hat{\boldsymbol{c}}_R^{(1)}\mu^2\right)t_{\mathrm{go}}^{(2)2} + \left(\hat{\boldsymbol{d}}_R^{(2)} + \hat{\boldsymbol{d}}_R^{(1)}\mu\right)t_{\mathrm{go}}^{(2)}
\end{aligned}
\tag{14}
$$

By examining Eq. (13), an auxiliary equality constraint on $\Delta\theta$ and an auxiliary equality constraint on $\Delta\psi$ can be added. This ensures that the selection of the $\Delta\theta$ function essentially does not affect the velocity increments in the *x* and *y* directions, and the selection of the $\Delta\psi$ function essentially does not affect the velocity increment in the *z* direction.

$$
\Theta\left(\Delta\theta\right) = 0 \tag{15}
$$

$$
\Theta\left(\Delta\psi\right) = 0 \tag{16}
$$

Substituting Eqs. (15) and (16) into Eq. (13), Eq. (13) can be simplified to:

$$
\mathbf{0} = \boldsymbol{v}_0 + \frac{T}{m_0}\begin{bmatrix} \cos\bar{\theta}\left[\chi - \frac{1}{2}\Theta\left(\Delta\theta^2\right)\right] \\ \sin\bar{\theta}\left[\chi - \frac{1}{2}\Theta\left(\Delta\theta^2\right)\right] \\ -\bar{\psi}\chi \end{bmatrix} + \left(\mu+1\right)t_{\mathrm{go}}^{(2)}\boldsymbol{g} - k\left\|\boldsymbol{v}_0\right\|t_{\mathrm{go}}^{(2)}\Pi\boldsymbol{v}_0 + \boldsymbol{\xi}
\tag{17}
$$

Similarly, the expressions for the terminal position constraints can be calculated through double integration. Doubly integrating Eq. (8) yields the terminal position constraint expressions:

$$\begin{bmatrix} 0 \\ r_{yt} \\ 0 \end{bmatrix} = \boldsymbol{r}_0 + \left[ (\mu+1) t_{\text{go}}^{(2)} - k \|\boldsymbol{v}_0\| t_{\text{go}}^{(2)2} \tilde{\Pi} \right] \boldsymbol{v}_0 + \frac{1}{2} (\mu+1)^2 t_{\text{go}}^{(2)2} \boldsymbol{g} + \tilde{\boldsymbol{\xi}}$$

$$+ \frac{T}{m_0} \begin{bmatrix} \cos\bar{\theta} \left[ \tilde{\chi} - \frac{1}{2} \tilde{\Theta}(\Delta\theta^2) \right] - \tilde{\Theta}(\Delta\theta) \sin\bar{\theta} \\ \sin\bar{\theta} \left[ \tilde{\chi} - \frac{1}{2} \tilde{\Theta}(\Delta\theta^2) \right] + \tilde{\Theta}(\Delta\theta) \cos\bar{\theta} \\ -\bar{\psi}\tilde{\chi} - \tilde{\Theta}(\Delta\psi) \end{bmatrix} \tag{18}$$

where

$$\begin{aligned} \tilde{\chi} &= \frac{T}{6 I_{\text{sp}} m_0} \left( n^2\mu^3 + 3n^2\mu^2 + 3n\mu + 1 \right) t_{\text{go}}^{(2)3} + \frac{1}{2} \left( n\mu^2 + 2n\mu + 1 \right) t_{\text{go}}^{(2)2} \\ \tilde{\Pi} &= \frac{1}{12(1+n\mu)^2} \left( 3 + 12\mu + 6\mu^2 + 12n\mu^2 + 8n\mu^3 + 4n^2\mu^3 + 3n^2\mu^4 \right) \\ \tilde{\boldsymbol{\xi}} &= \left( \frac{1}{20} \hat{\boldsymbol{a}}_R^{(1)} \mu^5 + \frac{1}{4} \hat{\boldsymbol{a}}_R^{(1)} \mu^4 + \frac{1}{20} \hat{\boldsymbol{a}}_R^{(2)} \right) t_{\text{go}}^{(2)5} + \left( \frac{1}{12} \hat{\boldsymbol{b}}_R^{(1)} \mu^4 + \frac{1}{3} \hat{\boldsymbol{b}}_R^{(1)} \mu^3 + \frac{1}{12} \hat{\boldsymbol{b}}_R^{(2)} \right) t_{\text{go}}^{(2)4} \\ &+ \left( \frac{1}{6} \hat{\boldsymbol{c}}_R^{(1)} \mu^3 + \frac{1}{2} \hat{\boldsymbol{c}}_R^{(1)} \mu^2 + \frac{1}{6} \hat{\boldsymbol{c}}_R^{(2)} \right) t_{\text{go}}^{(2)3} + \left( \frac{1}{2} \hat{\boldsymbol{d}}_R^{(1)} \mu^2 + \hat{\boldsymbol{d}}_R^{(1)} \mu + \frac{1}{2} \hat{\boldsymbol{d}}_R^{(2)} \right) t_{\text{go}}^{(2)2} \end{aligned} \tag{19}$$

In addition, the rocket landing must satisfy the terminal attitude verticality constraints:

$$\bar{\theta} + \Delta\theta \Big|_{\tau=(1+\mu)t_{\text{go}}^{(2)}} = \frac{\pi}{2} \tag{20}$$

$$\bar{\psi} + \Delta\psi \Big|_{\tau=(1+\mu)t_{\text{go}}^{(2)}} = 0 \tag{21}$$

Subsequently, based on constraints (15)~(21) and the objective function, this paper will derive the specific implementation method of the guidance algorithm.

## V. Specific Implementation Methodology

### A. Calculation of Key Parameters

The key parameters refer to the principal components of the pitch angle $\bar{\theta}$, the yaw angle $\bar{\psi}$, and the terminal time $t_{\text{go}}^{(2)}$. By rearranging the first two rows of Eq. (17) and eliminating $\bar{\theta}$ through squaring and summing, the equation can be simplified to:

$$\begin{aligned} &\frac{(n\mu+1)^4 T^4}{4 I_{\text{sp}}^2 m_0^4} t_{\text{go}}^{(2)4} + \frac{(n\mu+1)^3 T^3}{I_{\text{sp}} m_0^3} t_{\text{go}}^{(2)3} \\ &+ \left\{ (n\mu+1)^2 \frac{T^2}{m_0^2} - \frac{(n\mu+1)^2 T^3}{2 I_{\text{sp}} m_0^3} \Theta(\Delta\theta^2) - \left[ g(\mu+1) + k\|\boldsymbol{v}_0\| \Pi v_{y0} \right]^2 - k^2 \|\boldsymbol{v}_0\|^2 \Pi^2 v_{x0}^2 \right\} t_{\text{go}}^{(2)2} \\ &+ \left\{ -\frac{(n\mu+1) T^2}{m_0^2} \Theta(\Delta\theta^2) + 2\left[ g(\mu+1) + k\|\boldsymbol{v}_0\| \Pi v_{y0} \right] (v_{y0} + \xi_y) + 2k \|\boldsymbol{v}_0\| \Pi v_{x0} (v_{x0} + \xi_x) \right\} t_{\text{go}}^{(2)} \\ &+ \frac{T^2}{4 m_0^2} \Theta(\Delta\theta^2)^2 - (v_{y0} + \xi_y)^2 - (v_{x0} + \xi_x)^2 = 0 \end{aligned} \tag{22}$$

The above equation is a quartic equation in terms of the remaining flight time $t_{\text{go}}^{(2)}$, for which an analytical solution can be directly obtained. The specific method is detailed in Appendix C. It is worth noting that $\Theta\left(\Delta\theta^2\right)$ can be set to 0 in the first iteration and calculated using Eq. (B14) in subsequent iterations. Furthermore, based on the *x*-component in Eq. (17), the expression for the principal component of the pitch angle can be derived as:

$$\bar{\theta} = \arccos\left\{\frac{m_0}{T}\left(k\|\boldsymbol{v}_0\|t_{\mathrm{go}}^{(2)}\Pi v_{x0} - v_{x0} - \xi_x\right)\left[\chi - \frac{1}{2}\Theta\left(\Delta\theta^2\right)\right]^{-1}\right\} \tag{23}$$

Similarly, based on the *z*-component, the expression for the principal component of the yaw angle can be written as:

$$\bar{\psi} = \frac{m_0}{T}\left(v_{z0} + \xi_z - k\|\boldsymbol{v}_0\|t_{\mathrm{go}}^{(2)}\Pi v_{z0}\right)\chi^{-1} \tag{24}$$

So far, the analytical expressions for the principal components $\bar{\theta}$, $\bar{\psi}$, and the terminal time $t_{\mathrm{go}}^{(2)}$ have been derived. The terminal velocity component constraints (17) are satisfied. The dependency relationships between the parameters are illustrated in the following figure.

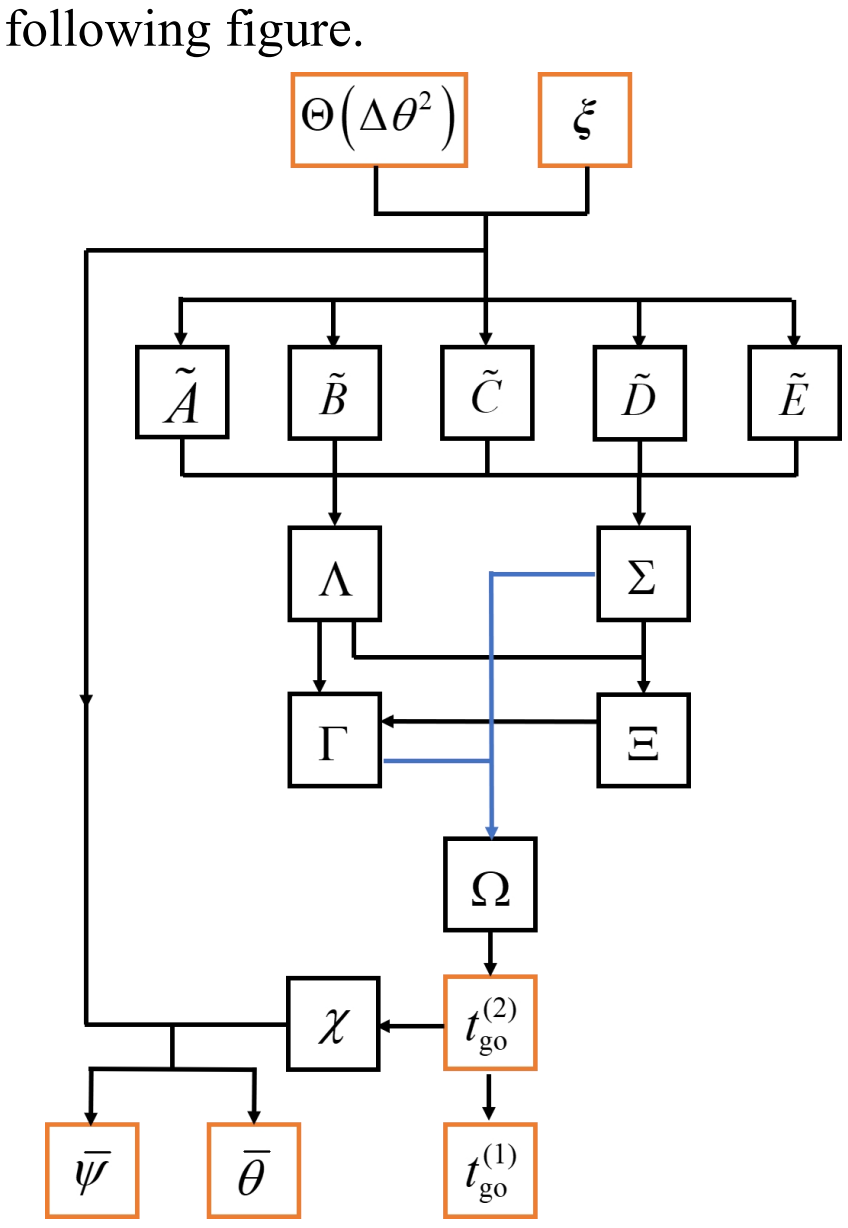

**Fig. 3 Dependency relationships among variables in the calculation of key parameters**

## B. Calculation of Pitch Angle Profile

Assuming the pitch angle increment $\Delta\theta$ is a piecewise cubic function (optimality proven), to ensure the rocket is vertically oriented at landing, the pitch angle profile can be defined as:

$$
\begin{aligned}
&\Delta\theta = \begin{cases} a^{(1)}\left(\tau - t_{\mathrm{go}}\right)^3 + b^{(1)}\left(\tau - t_{\mathrm{go}}\right)^2 + c^{(1)}\left(\tau - t_{\mathrm{go}}\right) + d^{(1)}, \tau \le \mu t_{\mathrm{go}}^{(2)} \\ a^{(2)}\left(\tau - t_{\mathrm{go}}\right)^3 + b^{(2)}\left(\tau - t_{\mathrm{go}}\right)^2 + c^{(2)}\left(\tau - t_{\mathrm{go}}\right) + d^{(2)}, \mu t_{\mathrm{go}}^{(2)} < \tau \le (\mu+1) t_{\mathrm{go}}^{(2)} \end{cases} \\
&a^{(1)} = n a^{(2)} \\
&b^{(1)} = n b^{(2)} \\
&c^{(1)} = c^{(2)} + (n-1) t_{\mathrm{go}}^{(2)} \left(2b^{(2)} - 3a^{(2)} t_{\mathrm{go}}^{(2)}\right) \\
&d^{(1)} = d^{(2)} + (n-1) t_{\mathrm{go}}^{(2)2} \left(b^{(2)} - 2a^{(2)} t_{\mathrm{go}}^{(2)}\right) \\
&d^{(2)} = \frac{\pi}{2} - \overline{\theta}
\end{aligned} \tag{25}
$$

This construction not only guarantees that the pitch angle is 90 degrees at the terminal moment but also ensures that the continuity, differentiability, and proportionality conditions listed in Theorem 1 are all satisfied.

Referring to Eq. (B5), constraint (15) can be equivalently written as:

$$A_1 a^{(2)} + B_1 b^{(2)} + C_1 c^{(2)} + D_1 d^{(2)} = 0 \tag{26}$$

where $A_1, B_1, C_1, D_1$ are parameters that can be calculated analytically (see Eq. (B6) in the Appendix for the specific formula). Next, the *x*-component in Eq. (18) is equivalently transformed into:

$$\tilde{\Theta}\left(\Delta\theta\right) = X \tag{27}$$

where

$$X = \frac{m_0}{T \sin\overline{\theta}} \left[ r_{x0} + v_{x0}(\mu+1) t_{\mathrm{go}}^{(2)} - k \left\| \boldsymbol{v}_0 \right\| t_{\mathrm{go}}^{(2)2} \tilde{\Pi} v_{x0} + \tilde{\xi}_x \right] + \frac{\cos\overline{\theta}}{\sin\overline{\theta}} \left[ \tilde{\chi} - \frac{1}{2} \tilde{\Theta}\left(\Delta\theta^2\right) \right] \tag{28}$$

Referring to Eq. (B11), constraint (27) can be equivalently written as:

$$A_2 a^{(2)} + B_2 b^{(2)} + C_2 c^{(2)} + D_2 d^{(2)} - X = 0 \tag{29}$$

where $A_2, B_2, C_2, D_2$ are parameters that can be calculated analytically (see Eq. (B12) in the Appendix for the specific formula). Next, by combining Eqs. (26) and (29), $b^{(2)}$ and $c^{(2)}$ can be expressed in terms of the parameter $a^{(2)}$. Solving this system of binary equations yields:

$$
\begin{aligned}
b^{(2)} &= B_a a^{(2)} + B_0 \\
c^{(2)} &= C_a a^{(2)} + C_0
\end{aligned} \tag{30}
$$

where

$$
\begin{aligned}
B_a &= \frac{A_1 C_2 - A_2 C_1}{B_2 C_1 - B_1 C_2}, B_0 = \frac{C_2 D_1 d^{(2)} - C_1 D_2 d^{(2)} + C_1 X}{B_2 C_1 - B_1 C_2} \\
C_a &= \frac{A_2 B_1 - A_1 B_2}{B_2 C_1 - B_1 C_2}, C_0 = \frac{B_1 D_2 d^{(2)} - B_1 X - B_2 D_1 d^{(2)}}{B_2 C_1 - B_1 C_2}
\end{aligned} \tag{31}
$$

At this point, considering that $d^{(2)}$ is a constant, the remaining parameters in the piecewise cubic function have all been expressed in terms of $a^{(2)}$. Next, the objective function is rearranged with $a^{(2)}$ as the independent variable. Adjusting $a^{(2)}$ to its minimum value ensures optimal smoothness of the pitch angle:

$$J\left(a^{(2)}\right)=\int_{\tau=0}^{(1+\mu)t_{\text{go}}^{(2)}}\left(\frac{\mathrm{d}\Delta\theta}{\mathrm{d}\tau}\right)^2\ \mathrm{d}\tau = Aa^{(2)2}+Ba^{(2)}+C \tag{32}$$

where

$$\begin{aligned} A &= \frac{9}{5}t_{\text{go}}^{(2)5}\kappa_5 - 3B_a t_{\text{go}}^{(2)4}\kappa_4 + t_{\text{go}}^{(2)3}\left(2C_a\kappa_{3,1}+\frac{4}{3}B_a^2\kappa_{3,2}\right) - 2B_aC_a t_{\text{go}}^{(2)2}\kappa_2 + C_a^2 t_{\text{go}}^{(2)}\kappa_1 \\ B &= -3B_0 t_{\text{go}}^{(2)4}\kappa_4 + t_{\text{go}}^{(2)3}\left(2C_0\kappa_{3,1}+\frac{8}{3}B_0B_a\kappa_{3,2}\right) - 2\left(B_aC_0+B_0C_a\right)t_{\text{go}}^{(2)2}\kappa_2 + 2C_0C_a t_{\text{go}}^{(2)}\kappa_1 \end{aligned} \tag{33}$$

and $\kappa_1 \sim \kappa_4$ are expansion coefficients greater than 1:

$$\begin{aligned} \kappa_1 &= 1+\mu \\ \kappa_2 &= 1+2\mu+n\mu^2 \\ \kappa_{3,1} &= 1+3\mu+3n\mu^2+n\mu^3 \\ \kappa_{3,2} &= 1+3\mu+3n\mu^2+n^2\mu^3 \\ \kappa_4 &= 1+4\mu+6n\mu^2+\frac{4}{3}n\mu^3+\frac{8}{3}n^2\mu^3+n^2\mu^4 \\ \kappa_5 &= 1+5\mu+10n\mu^2+\frac{10n\mu^3}{3}+\frac{20n^2\mu^3}{3}+5n^2\mu^4+n^2\mu^5 \end{aligned} \tag{34}$$

It can be observed that the objective function (32) is a quadratic function of $a^{(2)}$. The extreme point (which is also the minimum point) of this function is:

$$a^{(2)} = -\frac{B}{2A} \tag{35}$$

If $a^{(2)}, b^{(2)}, c^{(2)}$ satisfies the analytical expressions (35) and (30) respectively, then constraint Eqs. (15), (20), and the *x*-component of Eq. (18) are satisfied, and the pitch angle smoothness is optimal. Finally, the current pitch angle guidance command can be calculated according to the following formula:

$$\theta = -a^{(1)}t_{\text{go}}^3 + b^{(1)}t_{\text{go}}^2 - c^{(1)}t_{\text{go}} + d^{(1)} + \overline{\theta} \tag{36}$$

The dependency relationships between the parameters are illustrated in the following figure.

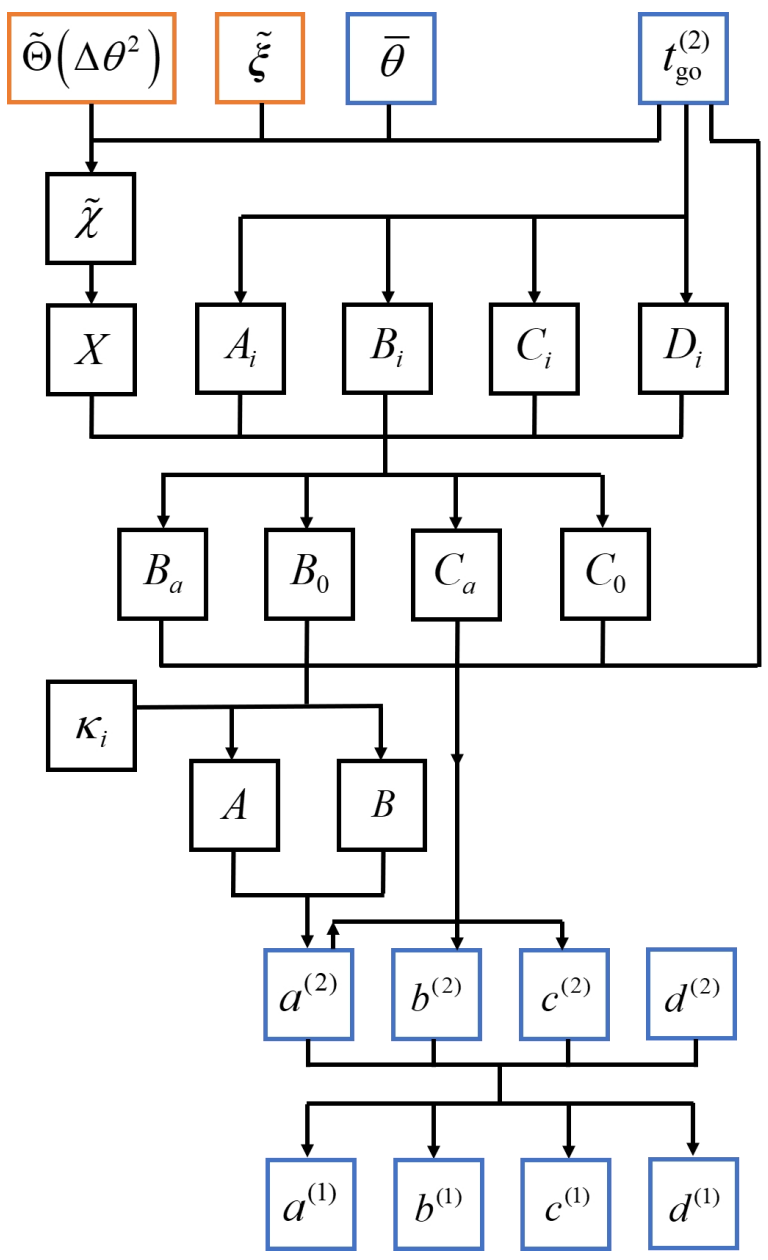

**Fig. 4 Dependency relationships among variables in the calculation of the pitch angle profile**

## C. Calculation of the Yaw Angle Profile

Similarly, it is assumed that the yaw angle increment is a piecewise cubic function satisfying the properties of Theorem 1:

$$\begin{aligned}
&\Delta\psi=\begin{cases} a'^{(1)}\left(\tau-t_{go}\right)^3+b'^{(1)}\left(\tau-t_{go}\right)^2+c'^{(1)}\left(\tau-t_{go}\right)+d'^{(1)}, \tau\le \mu t_{go}^{(2)} \\ a'^{(2)}\left(\tau-t_{go}\right)^3+b'^{(2)}\left(\tau-t_{go}\right)^2+c'^{(2)}\left(\tau-t_{go}\right)+d'^{(2)}, \mu t_{go}^{(2)}<\tau\le\left(\mu+1\right)t_{go}^{(2)} \end{cases}\\
&a'^{(1)}=na'^{(2)}\\
&b'^{(1)}=nb'^{(2)}\\
&c'^{(1)}=c'^{(2)}+\left(n-1\right)t_{go}^{(2)}\left(2b'^{(2)}-3a'^{(2)}t_{go}^{(2)}\right)\\
&d'^{(1)}=d'^{(2)}+\left(n-1\right)t_{go}^{(2)2}\left(b'^{(2)}-2a'^{(2)}t_{go}^{(2)}\right)\\
&d'^{(2)}=-\bar{\psi}
\end{aligned} \tag{37}$$

Due to the high similarity between the yaw and pitch channels, the detailed derivation is omitted herein, and the results are directly presented as follows:

$$\begin{aligned}
a'^{(2)}&=-\frac{B'}{2A'}\\
b'^{(2)}&=B'_a a'^{(2)}+B'_0\\
c'^{(2)}&=C'_a a'^{(2)}+C'_0
\end{aligned} \tag{38}$$

where

$$
\begin{aligned}
A' &= \frac{9}{5} t_{\mathrm{go}}^{(2)5} \kappa_5 - 3B_a' t_{\mathrm{go}}^{(2)4} \kappa_4 + t_{\mathrm{go}}^{(2)3} \left( 2C_a' \kappa_{3,1} + \frac{4}{3} B_a'^2 \kappa_{3,2} \right) - 2B_a' C_a' t_{\mathrm{go}}^{(2)2} \kappa_2 + C_a'^2 t_{\mathrm{go}}^{(2)} \kappa_1 \\
B' &= -3B_0' t_{\mathrm{go}}^{(2)4} \kappa_4 + t_{\mathrm{go}}^{(2)3} \left( 2C_0' \kappa_{3,1} + \frac{8}{3} B_0' B_a' \kappa_{3,2} \right) - 2 \left( B_a' C_0' + B_0' C_a' \right) t_{\mathrm{go}}^{(2)2} \kappa_2 + 2C_0' C_a' t_{\mathrm{go}}^{(2)} \kappa_1 \\
B_a' &= \frac{A_1 C_2 - A_2 C_1}{B_2 C_1 - B_1 C_2}, B_0' = \frac{C_2 D_1 d'^{(2)} - C_1 D_2 d'^{(2)} + C_1 Z}{B_2 C_1 - B_1 C_2} \\
C_a' &= \frac{A_2 B_1 - A_1 B_2}{B_2 C_1 - B_1 C_2}, C_0' = \frac{B_1 D_2 d'^{(2)} - B_1 Z - B_2 D_1 d'^{(2)}}{B_2 C_1 - B_1 C_2} \\
Z &= \frac{m_0}{T} \left[ r_{z0} + v_{z0} \left( \mu + 1 \right) t_{\mathrm{go}}^{(2)} - k \left\| \boldsymbol{v}_0 \right\| t_{\mathrm{go}}^{(2)2} \tilde{\Pi} v_{z0} + \tilde{\xi}_z \right] - \bar{\psi} \tilde{\chi}
\end{aligned}
\tag{39}
$$

If $a'^{(2)}, b'^{(2)}, c'^{(2)}$ satisfies the analytical expression (38), then constraint Eqs. (16), (21), and the $z$-component of Eq. (18) are satisfied, and the yaw angle smoothness is optimal. Finally, the current yaw angle guidance command can be calculated according to the following formula:

$$
\psi = -a'^{(1)} t_{\mathrm{go}}^3 + b'^{(1)} t_{\mathrm{go}}^2 - c'^{(1)} t_{\mathrm{go}} + d^{(1)} + \bar{\psi} \tag{40}
$$

The dependency relationships between the parameters are illustrated in the following figure.

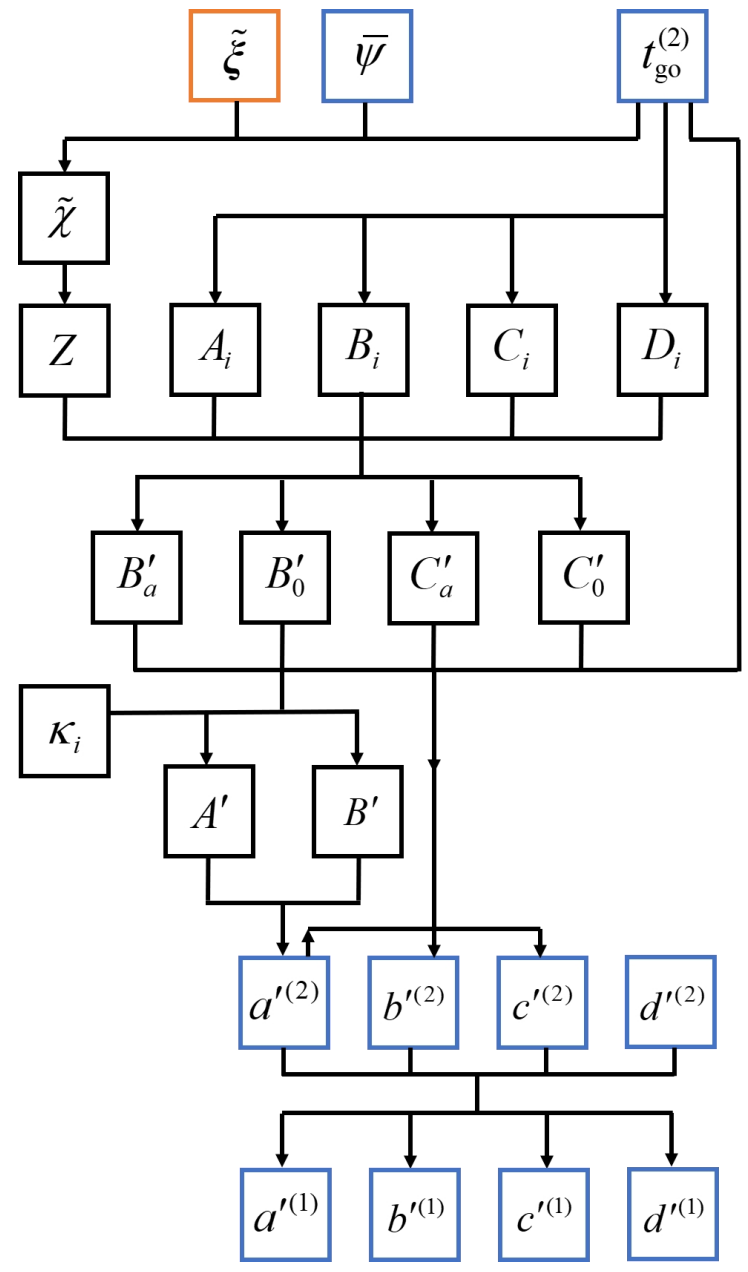

**Fig. 5 Dependency relationships among variables in the calculation of the yaw angle profile**

### D. Calculation of the Time Allocation Coefficient for the High-Thrust Phase

Thus far, Section A ensures the satisfaction of the three-axis terminal velocity constraints. Section B ensures the satisfaction of the $\Theta(\Delta\theta)$ constraints, the $x$-direction terminal position constraint, and the terminal pitch angle constraint. Section C ensures the satisfaction of the $\Theta(\Delta\psi)$ constraints, the $z$-direction terminal position constraint, and the terminal yaw angle constraint. Collectively, these satisfy the solution to the QLP. To ensure that a "fixed-point landing" can be achieved after guidance activation, the $y$-direction terminal position constraint must also be imposed. In other words, it is necessary to modify the time allocation coefficient $\mu$ or the thrust magnitude $T$ such that the rocket height decreases to 0 when the velocity decreases to 0. For the high-thrust phase, this paper chooses to satisfy the height constraint by modifying $\mu$, which completely avoids the need for thrust magnitude adjustment during this flight phase.

The $y$-direction terminal position constraint can be written as:

$$\begin{aligned} r_{yt} &= r_{y0} + v_{y0}\left(\mu+1\right)t_{go}^{(2)} - k\left\|\boldsymbol{v}_0\right\|t_{go}^{(2)2}\tilde{\Pi}v_{y0} - \frac{1}{2}g\left(\mu+1\right)^2 t_{go}^{(2)2} + \tilde{\xi}_y \\ &+ \frac{T}{m_0}\left\{\sin\bar{\theta}\left[\tilde{\chi} - \frac{1}{2}\tilde{\Theta}\left(\Delta\theta^2\right)\right] + \tilde{\Theta}\left(\Delta\theta\right)\cos\bar{\theta}\right\} \end{aligned} \tag{41}$$

Neglecting aerodynamic forces and mass variation, the approximate partial derivative of $r_{yt}$ with respect to $\mu$ can be found as:

$$\begin{aligned} \frac{\partial r_{yt}}{\partial \mu} &\approx v_{y0}t_{go}^{(2)} - g\left(\mu+1\right)t_{go}^{(2)2} + \frac{T}{m_0}\sin\bar{\theta}n\left(\mu+1\right)t_{go}^{(2)2} \\ &= v_{y0}t_{go}^{(2)} + \left(\frac{nT}{m_0}\sin\bar{\theta} - g\right)\left(\mu+1\right)t_{go}^{(2)2} \\ &= t_{go}^{(2)}\left[v_{y0} + \left(\frac{nT}{m_0}\sin\bar{\theta} - g\right)t_{go}\right] > 0 \end{aligned} \tag{42}$$

The second factor in Eq. (42) approximately characterizes the vertical velocity component of the rocket after applying the high-thrust mode for $t_{go}$, which is generally a positive value. This indicates that increasing $\mu$ will result in an increase in the terminal height $r_{yt}$. The correction formula based on Newton's iteration can be adopted as follows:

$$\mu \leftarrow \mu - \frac{\alpha}{t_{go}^{(2)}}\left[v_{y0} + \left(\frac{nT}{m_0}\sin\bar{\theta} - g\right)t_{go}\right]^{-1} r_{yt} \tag{43}$$

where $\alpha$ is a step size less than 1. After multiple iterative corrections using formula (43), the converged $\mu$ ensures that the rocket terminal height $r_{yt} = 0$.

**E. Calculation of the Thrust Magnitude for the Low-Thrust Phase**

After entering the low-thrust phase, the value of $\mu$ is fixed at 0, making it impossible to achieve a fixed-point landing by adjusting $\mu$. To ensure a fixed-point landing during the low-thrust phase, the thrust magnitude must be modified. By directly rearranging Eq. (41), the correction formula for $T$ can be obtained:

$$T \leftarrow m_0 \frac{\frac{1}{2}gt_{go}^{(2)2} + k\left\|\boldsymbol{v}_0\right\|t_{go}^{(2)2}\tilde{\Pi}v_{y0} - \tilde{\xi}_y - r_{y0} - v_{y0}t_{go}^{(2)}}{\sin\bar{\theta}\left[\tilde{\chi} - \frac{1}{2}\tilde{\Theta}\left(\Delta\theta^2\right)\right] + \tilde{\Theta}\left(\Delta\theta\right)\cos\bar{\theta}} \tag{44}$$

where $\tilde{\Theta}(\Delta\theta)$ can be obtained from Eq. (B11), $\tilde{\Theta}(\Delta\theta^2)$ from Eq. (B15), and $\bar{\theta}$ from Eq. (23).

### F. Aerodynamic Correction Method

This section is primarily used to calculate the undetermined parameters of $\Delta\boldsymbol{a}_R(\tau)$ in Eq. (7) to estimate aerodynamic acceleration. It can be seen that $\Delta\boldsymbol{a}_R(\tau)$ has a total of 8 undetermined parameters. Therefore, 8 linear equations can be constructed to determine these parameters. First, 4 key time points are selected: the initial moment, the moment immediately before the thrust switching, the moment immediately after the thrust switching, and the terminal moment. An auxiliary function is constructed as:

$$\boldsymbol{\eta}(\tau)=\tilde{\boldsymbol{a}}_R\left[\theta(\tau),\boldsymbol{v}(\tau)\right]-\bar{\boldsymbol{a}}_R(\tau) \tag{45}$$

Then, according to Eq. (7), a total of 8 constraint equations can be written as follows:

$$\begin{aligned}&\Delta\boldsymbol{a}_R(\tau)=\boldsymbol{\eta}(\tau)\\&\Delta\dot{\boldsymbol{a}}_R(\tau)=\dot{\boldsymbol{\eta}}(\tau)\\&\tau=0,\mu t_{\text{go}}^{(2)-},\mu t_{\text{go}}^{(2)+},(\mu+1)t_{\text{go}}^{(2)}\end{aligned} \tag{46}$$

If Eq. (46) is satisfied, it ensures that the aerodynamic acceleration and its rate of change in the simplified dynamic model are consistent with those in the precise dynamic model at the four key time points, achieving a relatively accurate estimation of the aerodynamic force. The $\dot{\boldsymbol{\eta}}(\tau)$ at each key moment can be obtained by differentiating Eq. (45) with respect to time:

$$\dot{\boldsymbol{\eta}}(\tau)=\dot{\tilde{\boldsymbol{a}}}_R\left[\theta(\tau),\boldsymbol{v}(\tau)\right]-\dot{\bar{\boldsymbol{a}}}_R(\tau) \tag{47}$$

where

$$\begin{gathered}\dot{\tilde{\boldsymbol{a}}}_R\left[\theta(\tau),\boldsymbol{v}(\tau)\right]=-k\frac{\boldsymbol{v}^T\dot{\boldsymbol{v}}}{\|\boldsymbol{v}\|}\begin{bmatrix}f_x(\theta,v_x,v_y)\\f_y(\theta,v_x,v_y)\\f_z(v_z)\end{bmatrix}-k\|\boldsymbol{v}\|\begin{bmatrix}f_x(\theta,\dot{v}_x,\dot{v}_y)+\left(1-\dfrac{C_y}{C_x}\right)\dot{\theta}\left(-v_x\sin2\theta+v_y\cos2\theta\right)\\f_y(\theta,\dot{v}_x,\dot{v}_y)+\left(1-\dfrac{C_y}{C_x}\right)\dot{\theta}\left(v_y\sin2\theta+v_x\cos2\theta\right)\\f_z(\dot{v}_z)\end{bmatrix}\\ \dot{\bar{\boldsymbol{a}}}_R(\tau)=\begin{cases}\left[1-\dfrac{n\tau}{(n\mu+1)t_{\text{go}}^{(2)}}\right]\dfrac{2nk\|\boldsymbol{v}_0\|}{(n\mu+1)t_{\text{go}}^{(2)}}\boldsymbol{v}_0, & \tau\le\mu t_{\text{go}}^{(2)}\\ \left[\dfrac{(\mu+1)t_{\text{go}}^{(2)}-\tau}{(n\mu+1)t_{\text{go}}^{(2)}}\right]\dfrac{2k\|\boldsymbol{v}_0\|}{(n\mu+1)t_{\text{go}}^{(2)}}\boldsymbol{v}_0, & \mu t_{\text{go}}^{(2)}<\tau\le(\mu+1)t_{\text{go}}^{(2)}\end{cases}\end{gathered} \tag{48}$$

The polynomial coefficients can then be obtained by solving Eq. (46). The aerodynamic correction polynomial coefficients for the high-thrust phase are:

$$
\begin{aligned}
&\tau_1 = 0, \tau_2 = t_{\mathrm{go}}^{(2)-}\mu \\
&\hat{\boldsymbol{a}}_R^{(1)} = 2\frac{\boldsymbol{\eta}(\tau_1)-\boldsymbol{\eta}(\tau_2)}{t_{\mathrm{go}}^{(2)2}\mu^2} + \frac{\dot{\boldsymbol{\eta}}(\tau_1)+\dot{\boldsymbol{\eta}}(\tau_2)}{t_{\mathrm{go}}^{(2)}\mu} \\
&\hat{\boldsymbol{b}}_R^{(1)} = 3\frac{-\boldsymbol{\eta}(\tau_1)+\boldsymbol{\eta}(\tau_2)}{t_{\mathrm{go}}^{(2)2}\mu^2} - \frac{2\dot{\boldsymbol{\eta}}(\tau_1)+\dot{\boldsymbol{\eta}}(\tau_2)}{t_{\mathrm{go}}^{(2)}\mu} \\
&\hat{\boldsymbol{c}}_R^{(1)} = \dot{\boldsymbol{\eta}}(\tau_1) \\
&\hat{\boldsymbol{d}}_R^{(1)} = \boldsymbol{\eta}(\tau_1)
\end{aligned} \tag{49}
$$

The aerodynamic correction polynomial coefficients for the low-thrust phase are:

$$
\begin{aligned}
&\tau_1 = t_{\mathrm{go}}^{(2)+}\mu, \tau_2 = t_{\mathrm{go}}^{(2)}(1+\mu) \\
&\hat{\boldsymbol{a}}_R^{(2)} = 2\frac{\boldsymbol{\eta}(\tau_1)-\boldsymbol{\eta}(\tau_2)}{t_{\mathrm{go}}^{(2)2}} + \frac{\dot{\boldsymbol{\eta}}(\tau_1)+\dot{\boldsymbol{\eta}}(\tau_2)}{t_{\mathrm{go}}^{(2)}} \\
&\hat{\boldsymbol{b}}_R^{(2)} = 3\frac{-\boldsymbol{\eta}(\tau_1)+\boldsymbol{\eta}(\tau_2)}{t_{\mathrm{go}}^{(2)2}} - \frac{2\dot{\boldsymbol{\eta}}(\tau_1)+\dot{\boldsymbol{\eta}}(\tau_2)}{t_{\mathrm{go}}^{(2)}} \\
&\hat{\boldsymbol{c}}_R^{(2)} = \dot{\boldsymbol{\eta}}(\tau_1) \\
&\hat{\boldsymbol{d}}_R^{(2)} = \boldsymbol{\eta}(\tau_1)
\end{aligned} \tag{50}
$$

Specifically, the aerodynamic correction polynomial coefficients can be calculated in sequence according to the following method, referred to as **Algorithm 1**:

**Algorithm 1:** Update Algorithm for Aerodynamic Correction Coefficients

Step1: Input $t_{\mathrm{go}}^{(2)}$ , $\mu$ , $\bar{\theta}$ , $\boldsymbol{v}_0$ and necessary parameters, initial guess of aerodynamic correction coefficients $\hat{\boldsymbol{a}}_R^{(i)} \sim \hat{\boldsymbol{d}}_R^{(i)}$.

For $\tau=0, \mu t_{\mathrm{go}}^{(2)-}, \mu t_{\mathrm{go}}^{(2)+}, (\mu+1)t_{\mathrm{go}}^{(2)}$

  Step2: Calculate $\theta(\tau)$ according to Eq. (25) and $\bar{\theta}$ .

  Step3: Calculate pitch angle rate $\dot{\theta}(\tau)$ according to the derivative of Eq. (25).

  If $\tau=0$

    Step4: $\boldsymbol{v}(\tau) = \boldsymbol{v}_0$

  Elseif $\tau=\mu t_{\mathrm{go}}^{(2)-}$ or $\tau=\mu t_{\mathrm{go}}^{(2)+}$

    Step5: Calculate $\boldsymbol{v}(\tau) = \boldsymbol{v}_{\mathrm{m}}$ according to Eqs. (11), (B8), and (B16).

  Elseif $\tau=(\mu+1)t_{\mathrm{go}}^{(2)}$

    Step6: $\boldsymbol{v}(\tau) = 0$

  End

  Step7: Calculate $\tilde{\boldsymbol{a}}_R[\theta(\tau), \boldsymbol{v}(\tau)]$ and $\bar{\boldsymbol{a}}_R(\tau)$ according to Eqs. (5) and (7), respectively.

  Step8: Calculate velocity rate $\dot{\boldsymbol{v}}(\tau)$ according to Eq. (8).

  Step9: Calculate $\dot{\tilde{\boldsymbol{a}}}_R[\theta(\tau), \boldsymbol{v}(\tau)]$ and $\dot{\bar{\boldsymbol{a}}}_R(\tau)$ according to Eq. (48).

  Step10: Calculate $\boldsymbol{\eta}(\tau)$ and $\dot{\boldsymbol{\eta}}(\tau)$ according to Eqs. (45) and (47), respectively.

End

Step11: Update aerodynamic correction polynomial coefficients for the high-thrust phase $\hat{\boldsymbol{a}}_R^{(1)} \sim \hat{\boldsymbol{d}}_R^{(1)}$ according to Eq. (49).

Step12: Update aerodynamic correction polynomial coefficients for the low-thrust phase $\hat{\boldsymbol{a}}_R^{(2)} \sim \hat{\boldsymbol{d}}_R^{(2)}$ according to Eq. (50).

## VI. Summary of the Algorithm

The following is a summary of the proposed algorithm. This paper integrates the pre-activation phase (flag = 0), the high-thrust phase (flag = 1), and the low-thrust phase (flag = 2) into a complete algorithmic framework, referred to as **Algorithm 2**.

**Algorithm 2:** Smooth Attitude Maneuvering with Dual-thrust Switching (SAM-DS)

Step1: Set initial integral variables $\Theta(\Delta\theta^2)=0,\tilde{\Theta}(\Delta\theta^2)=0$, time allocation ratio $\mu$, nominal thrust magnitude *T*, initial aerodynamic correction coefficients $\hat{\boldsymbol{a}}_R^{(i)}=\boldsymbol{0},\hat{\boldsymbol{b}}_R^{(i)}=\boldsymbol{0},\hat{\boldsymbol{c}}_R^{(i)}=\boldsymbol{0},\hat{\boldsymbol{d}}_R^{(i)}=\boldsymbol{0}$, set the initial phase $\text{flag}=0$ and configure other necessary parameters.

For $iter=1:iter_{\max}$

Step2: Acquire the current state $\boldsymbol{r}_0$ $\boldsymbol{v}_0$ $m_0$

For $i=1:i_{\max}$

Step3: Calculate the coefficients of the remaining flight time equation $\tilde{A}$ $\tilde{B}$ $\tilde{C}$ $\tilde{D}$ $\tilde{E}$ according to Eq. (C1).

Step4: Calculate the remaining flight time $t_{\text{go}}^{(2)}$ according to Eq. (C2).

Step5: Calculate the coefficients of the constraint equations $A_i$ $B_i$ $C_i$ $D_i$ according to Eqs. (B6) and (B12).

| Pitch Channel | Yaw Channel |
|---|---|
| Step6: Calculate the pitch angle principal component $\bar{\theta}$ according to Eq. (23). | Calculate the yaw angle principal component $\bar{\psi}$ according to Eq. (24). |
| Step7: Calculate the pitch channel correlation parameter $B_a$ $B_0$ $C_a$ $C_0$ according to Eq. (31). | Calculate the yaw channel correlation parameter $B_a'$ $B_0'$ $C_a'$ $C_0'$ according to Eq. (39). |
| If $\text{flag}=0$ | |
| Step8: Calculate the pitch channel objective function coefficients $A$ $B$ according to Eq. (33). | Calculate the yaw channel objective function coefficients $A'$ $B'$ according to Eq. (39). |
| Step9: Calculate the pitch angle profile parameters $a^{(2)}$ $b^{(2)}$ $c^{(2)}$ according to Eqs. (35) and (30). | Calculate the yaw angle profile parameters $a'^{(2)}$ $b'^{(2)}$ $c'^{(2)}$ according to Eq. (38). |
| Else | |
| Step10: Fix $a^{(2)}$, calculate the pitch angle profile parameters $b^{(2)}$ $c^{(2)}$ according to Eqs. (35) and (30). | Fix $a'^{(2)}$, calculate the yaw angle profile parameters $b'^{(2)}$ $c'^{(2)}$ according to Eq. (38). |
| End | |
| Step11: Calculate the pitch angle profile $a^{(1)}$ $b^{(1)}$ $c^{(1)}$ $d^{(1)}$ according to Eq. (25). | Calculate the yaw angle profile $a'^{(1)}$ $b'^{(1)}$ $c'^{(1)}$ $d'^{(1)}$ according to Eq. (37). |
| Step12: Calculate the current pitch angle command $\theta$ according to Eq. (36). | Calculate the current yaw angle command $\psi$ according to Eq. (40). |

Step13: Update the integral variables $\Theta(\Delta\theta^2),\tilde{\Theta}(\Delta\theta^2)$ according to Eqs. (B14) and (B15).

Step14: Update the aerodynamic correction polynomial coefficients $\hat{\boldsymbol{a}}_R^{(i)},\hat{\boldsymbol{b}}_R^{(i)},\hat{\boldsymbol{c}}_R^{(i)},\hat{\boldsymbol{d}}_R^{(i)}$ according to **Algorithm 1**.

Step15: Calculate the terminal flight height $r_{yt}$ according to Eq. (41).

Step16: Calculate the variation of the integral variables between the current and previous iterations $\delta\Theta(\Delta\theta^2)$ $\delta\tilde{\Theta}(\Delta\theta^2)$

If $\text{flag}=1$

Step17: Update the time allocation coefficient $\mu$ according to Eq. (43).

Step18: Set the guidance thrust: $\text{Thrust} \leftarrow nT$ .
Elseif $\text{flag} = 2$
Step19: Update the thrust magnitude $T$ according to Eq. (44).
Step20: Set the guidance thrust: $\text{Thrust} \leftarrow T$ .
Else
Step21: Set the guidance thrust: $\text{Thrust} \leftarrow 0$ .
End
If $\left|\delta\Theta\left(\Delta\theta^2\right)\right| \le \varepsilon_\Theta$ $\left|\delta\tilde{\Theta}\left(\Delta\theta^2\right)\right| \le \varepsilon_{\tilde{\Theta}}$
Break;
End
End
Step22: Output the guidance commands $(\theta, \psi, \text{Thrust})$ .
If $\text{flag} = 0$ and $r_{\text{yt}} \le \varepsilon_{yt}$
Step23: $\text{flag} \leftarrow 1$
Elseif $\text{flag} = 1$ and $\mu \le \varepsilon_\mu$
Step24: $\mu \leftarrow \varepsilon_\mu$ , $\text{flag} \leftarrow 2$
Elseif $\text{flag} = 2$ and $t_{\text{go}}^{(2)} \le \varepsilon_t$
Break;
End
Step25: Wait for one guidance cycle.
End
Step26: Complete fixed-point landing. Algorithm terminates.

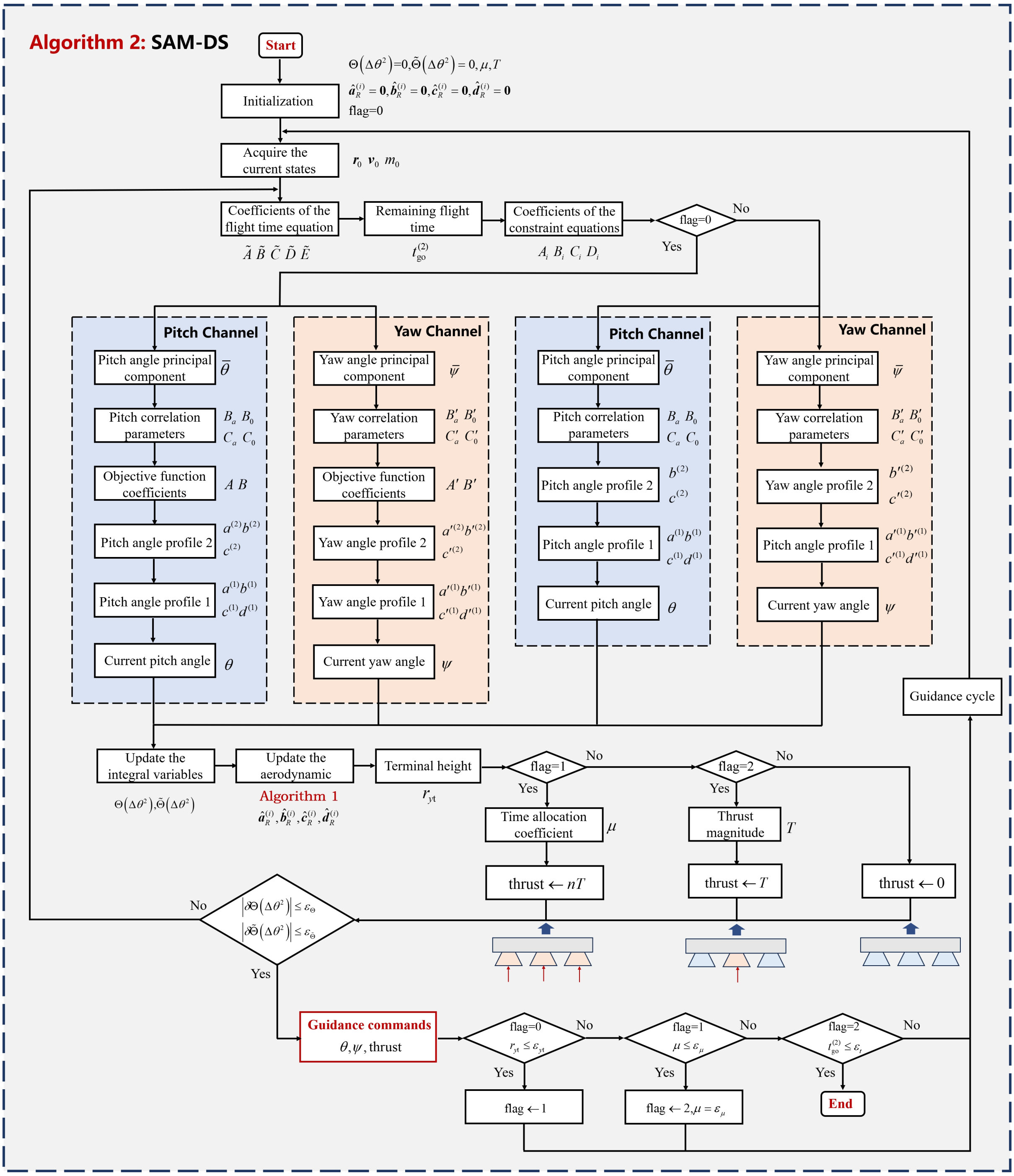

**Fig. 6 Algorithm Flowchart of SAM-DS**

The main differences under different flag values are summarized as follows:

(1) When flag=0 or 1, the thrust magnitude $T$ is set to a given constant. To provide a certain control margin, this constant value is typically set slightly below the maximum thrust of the engine during actual simulations. In contrast, when flag=2, to accurately satisfy the terminal height constraint, the thrust magnitude $T$ requires online correction (Step 18). Therefore, $T$ is treated as an undetermined constant to be solved in each guidance cycle. In practical applications, due to dynamic model errors, this thrust value will fluctuate slightly around the nominal thrust, which is a normal adjustment process of the SAM-DS algorithm.

(2) When flag=0, the profile parameters $a^{(2)}$ and $a'^{(2)}$ are determined by solving the extremum problem of a quadratic function, completing the initial trajectory planning. When flag switches from 0 to 1, the parameters $a^{(2)}$ and $a'^{(2)}$ are immediately frozen and treated as fixed values. Subsequently, the core task of the algorithm shifts to online updating of the remaining profile parameters. This approach maximizes the stability of the attitude angle profile and avoids unnecessary attitude fluctuations caused by adjusting multiple parameters simultaneously.

The complete flowchart of the guidance algorithm is shown in Fig. 6. It should be noted that Fig. 6 represents the full implementation of the logic outlined in Fig. 2.

## VII. Numerical Simulation

### A. Optimality and Real-Time Performance Analysis

To rigorously evaluate the computational effectiveness of the proposed method, two optimal control problems are defined as comparison benchmarks, named the "Approximate Problem" and the "Precise Problem," respectively. The "Approximate Problem" corresponds to Problem 1 defined previously, whose dynamic model is based on simplified assumptions. The "Precise Problem," however, is constructed by replacing the dynamic equations in Problem 1 with the Eq. (3). To maintain rigorous presentation, the complete mathematical description of the "Precise Problem" is defined as follows, referred to as **Problem 2**:

**Problem 2**:

$$\underset{u1,u2,u3,t_{\mathrm{go}}^{(2)}}{\text{Minimize}} \ J = \int_{\tau=0}^{(1+\mu)t_{\mathrm{go}}^{(2)}} \frac{1}{2}\left(u_1^2 + u_2^2\right) \mathrm{d}\,\tau$$

$$\dot{\boldsymbol{r}} = \boldsymbol{v},\ \dot{\boldsymbol{v}} = \boldsymbol{a}_T\left(T,\theta,\psi,\tau\right) + \boldsymbol{a}_R\left(\theta,\psi,\boldsymbol{v},\tau\right) + \boldsymbol{g},\ \dot{T} = 0$$

$$\dot{\theta} = u_1,\ \dot{\psi} = u_2,\ \dot{T} = u_3,\ u_3 = 0$$

$$\left[\boldsymbol{r}^T, \boldsymbol{v}^T\right]^T \Big|_{\tau=0} = \left[\boldsymbol{r}_0^T, \boldsymbol{v}_0^T\right]^T, \left[\boldsymbol{r}^T, \boldsymbol{v}^T\right]^T \Big|_{\tau=(1+\mu)t_{\mathrm{go}}^{(2)}} = \boldsymbol{0}_{6\times 1}$$

$$\theta\Big|_{\tau=(1+\mu)t_{\mathrm{go}}^{(2)}} = \frac{\pi}{2}, \psi\Big|_{\tau=(1+\mu)t_{\mathrm{go}}^{(2)}} = 0$$

Both Problem 1 and Problem 2 can be directly solved using the Sequential Convex Programming (SCP) method. Given that the SCP method is relatively mature, detailed elaboration is omitted here. The following will focus on comparing the differences among the proposed analytical method, the optimal solution of Problem 1, and the optimal solution of Problem 2. The parameter settings adopted for the simulation are listed in Table 1.

**Table 1 Parameters for Open-Loop Simulation**

| Item | Symbol | Value | Item | Symbol | Value |
|---|---|---|---|---|---|
| Nominal Thrust | $T$ | 700 kN | Initial Position | $\boldsymbol{r}_0$ | $[-400,1277,15]^T$ m |
| Specific Impulse | $I_{sp}$ | 3038 m/s | Initial Velocity | $\boldsymbol{v}_0$ | $[-20,-150,15]^T$ m/s |
| Atmospheric Density | $\rho$ | 1.225 kg/m$^3$ | Initial Mass | $m_0$ | 60 t |
| Reference Area | $S$ | 10 m$^2$ | Gravitational Acceleration | $g$ | 9.8 m/s$^2$ |
| $x$-axis Aerodynamic Coefficient | $C_x$ | 1 | Integration Step Size | $\Delta t$ | $10^{-4}$ s |
| $y$-axis Aerodynamic Coefficient | $C_y$ | 7 | Convex Optimization Discretization Steps | $N$ | 50 |
| $z$-axis Aerodynamic Coefficient | $C_z$ | 7 | Integral Parameter 1 Threshold | $\varepsilon_{\Theta}$ | $10^{-6}$ s |
| Thrust Ratio | $n$ | 3 | Integral Parameter 2 Threshold | $\varepsilon_{\tilde{\Theta}}$ | $10^{-6}$ s |
| Time Allocation Ratio | $\mu$ | 0.5 | | | |

The remaining flight times $t_{go}^{(2)}=14.23$ s and $t_{go}^{(1)}=7.11$ s are calculated according to Eq. (C2). The principal components of the pitch and yaw angles, calculated using Eqs. (23) and (24), are 79.4014 deg and 1.7201 deg, respectively. Subsequently, the pitch and yaw angle profiles are obtained using Eqs. (30), (35), and (38). The specific results are as follows:

For the pitch channel:

$$\begin{aligned} &a^{(1)} = -8.467\times10^{-4}, b^{(1)} = -0.03916, c^{(1)} = -0.5174, d^{(1)} = -1.848 \\ &a^{(2)} = -2.822\times10^{-4}, b^{(2)} = -0.01305, c^{(2)} = -0.1172, d^{(2)} = 0.1849 \end{aligned} \tag{51}$$

For the yaw channel:

$$\begin{aligned} &a'^{(1)} = 1.070\times10^{-4}, b'^{(1)} = 4.923\times10^{-3}, c'^{(1)} = 0.06422, d'^{(1)} = 0.2233 \\ &a'^{(2)} = 3.568\times10^{-5}, b'^{(2)} = 1.641\times10^{-3}, c'^{(2)} = 0.01416, d'^{(2)} = -0.03002 \end{aligned} \tag{52}$$

The curves of the pitch and yaw angles with time, as well as the trajectory, are depicted in Fig. 7.

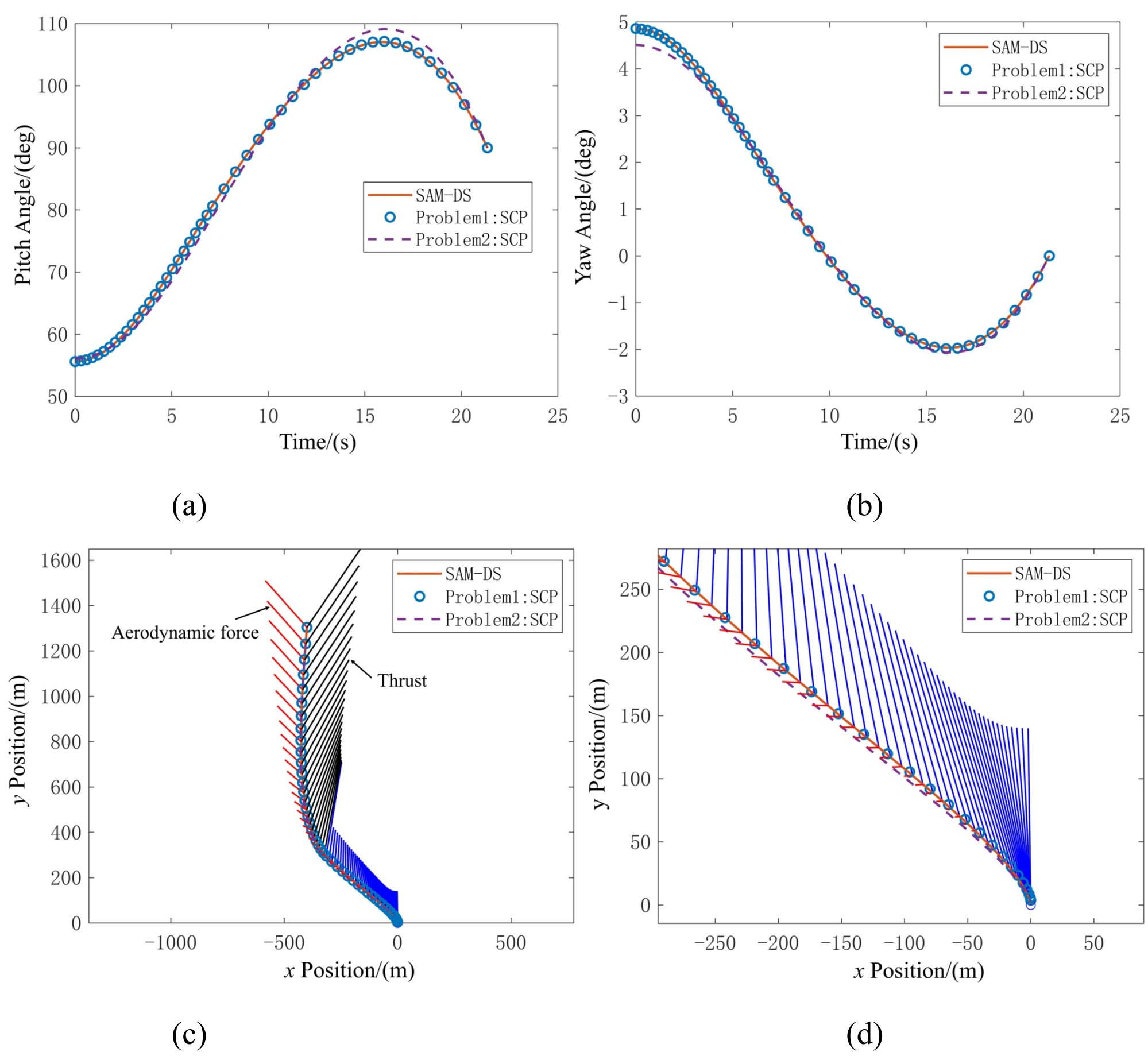

**Fig. 7 Time-history curves of control variables in open-loop simulation. (a) Pitch angle vs. time. (b) Yaw angle vs. time. (c) Trajectory curve. (d) Locally magnified trajectory curve.**

It can be observed from Fig. 7 that the proposed analytical solution is highly consistent with the optimal solution of the approximate model, indicating that the preconditions of Theorem 1 are basically satisfied. Specifically:

$$\frac{(n-1)T}{I_{\text{sp}}m_0}\mu t_{\text{go}}^{(2)}+\frac{T}{I_{\text{sp}}m_0}\tau\le 0.1366,\left|\frac{1}{2}\cos\bar{\theta}\Delta\theta^2\right|\le 0.0343,\left|\frac{1}{2}\sin\bar{\theta}\Delta\theta^2\right|\le 0.1834 \tag{53}$$

It can also be observed from Fig. 7 that the analytical solution is generally consistent with the optimal solution of the precise model, although certain discrepancies exist. These differences mainly stem from the simplified treatment of aerodynamic forces. In this paper, the aerodynamic forces are approximated using Algorithm 1. Due to the limited degrees of freedom of the polynomial, there is a slight deviation from the actual aerodynamic forces, leading to the differences between the analytical and precise solutions. The time-history curves of the aerodynamic acceleration predicted by the SAM-DS algorithm and the actual aerodynamic acceleration can be plotted, as shown in Fig. 8.

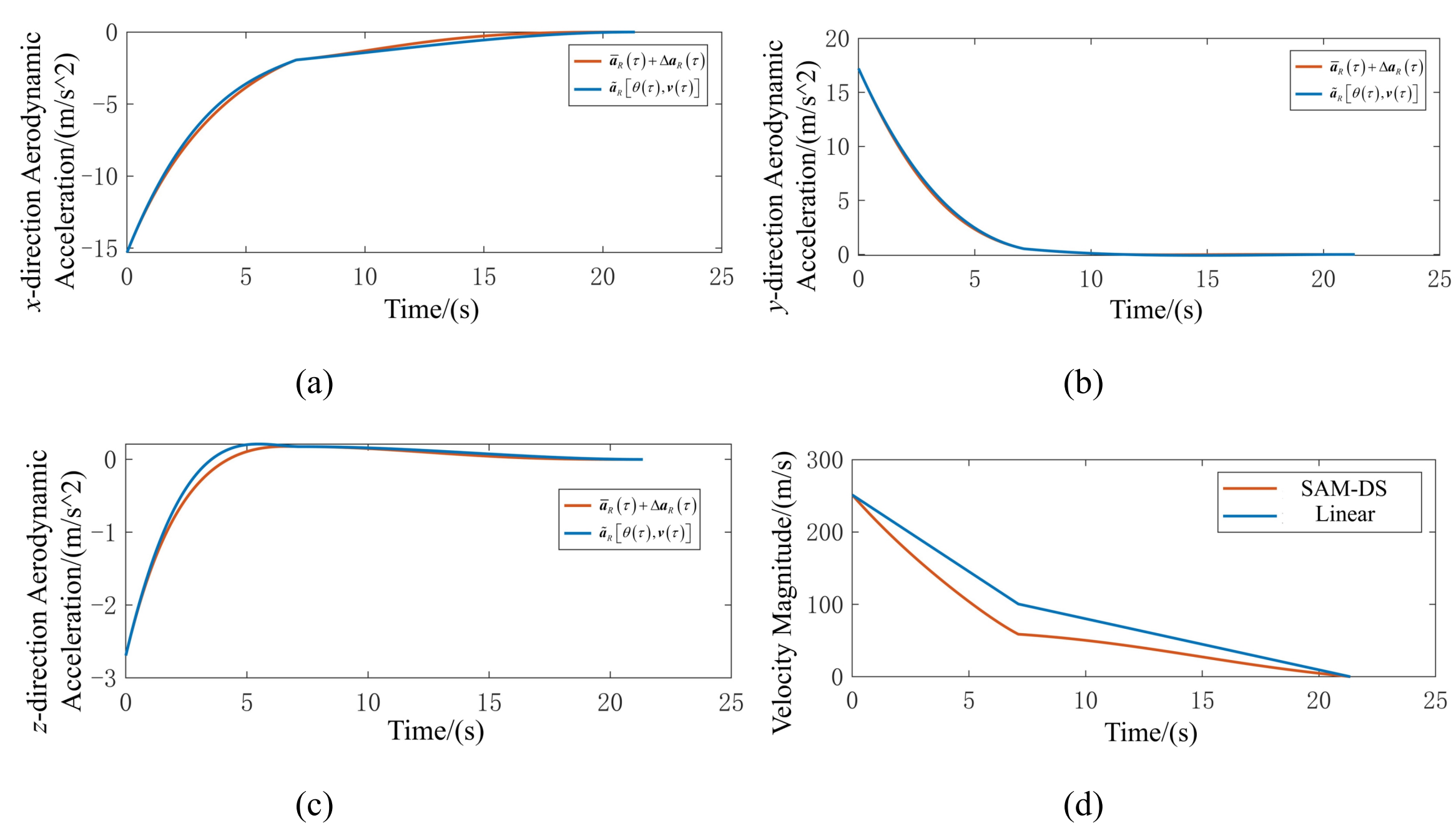

**Fig. 8 Time-history curves of aerodynamic acceleration. (a) *x*-component. (b) *y*-component. (c) *z*-component. (d) Piecewise linear assumption.**

As seen in Fig. 8 (d), due to the complexity of aerodynamic forces, the linear velocity assumption model exhibits significant discrepancies. If only the $\bar{\boldsymbol{a}}_R(\tau)$ term is considered while neglecting the $\Delta\boldsymbol{a}_R(\tau)$ term during the guidance process, large deviations will be introduced. This illustrates the necessity of estimating $\Delta\boldsymbol{a}_R(\tau)$. Observing Fig. 8 (a)~(c), it can be seen that after interpolation estimation using Algorithm 1, the three-axis aerodynamic acceleration $\bar{\boldsymbol{a}}_R(\tau)+\Delta\boldsymbol{a}_R(\tau)$ is basically consistent with $\tilde{\boldsymbol{a}}_R\left[\theta(\tau),\boldsymbol{v}(\tau)\right]$, reflecting the effectiveness of Algorithm 1.

Finally, a brief analysis of the computational real-time performance is provided. In this example, the simulation is conducted in the MATLAB environment on a desktop computer equipped with an Intel Core i7-11700K processor (3.60 GHz). The SAM-DS algorithm converges after 25 iterations, with a computational time of 0.56 ms, demonstrating excellent real-time performance. Compiling the algorithm in C or implementing a warm-start strategy will further improve the computational efficiency.

### B. Monte Carlo Simulation

Next, the robustness of the proposed algorithm is verified through systematic Monte Carlo simulations. The settings for random deviation parameters are shown in Table 4. Each parameter follows a uniform distribution, and the number of simulation samples is 1000. It is particularly noted that to fully investigate the algorithm's adaptability to large deviations in aerodynamic parameters under high lift-to-drag ratio conditions, the nominal lift-to-drag ratio is set to 7 in Table 4, and the aerodynamic parameter error is set to 50%.

**Table 4 Parameter Settings for Monte Carlo Simulation**

| Item | Symbol | Value | Item | Symbol | Value |
|---|---|---|---|---|---|
| Nominal Thrust | $T$ | 700 kN | Initial Position | $\boldsymbol{r}_0$ | $[0,2500,0]^T \pm [400,0,100]^T$ m |
| Specific Impulse | $I_{sp}$ | $3038 \pm 150$ m/s | Initial Velocity | $\boldsymbol{v}_0$ | $[0,-250,0]^T \pm [20,20,15]^T$ m |
| Atmospheric Density | $\rho$ | $1.225 \pm 0.1$ kg/m$^3$ | Initial Mass | $m_0$ | $60 \pm 4$ t |
| Reference Area | $S$ | $10 \pm 1$ m$^2$ | Gravitational Acceleration | $g$ | 9.8 m/s$^2$ |
| $x$-axis Aerodynamic Coefficient | $C_x$ | $1 \pm 0.5$ | Integration Step Size | $\Delta t$ | $10^{-4}$ s |
| $y$-axis Aerodynamic Coefficient | $C_y$ | $7 \pm 3.5$ | Integral Parameter 1 Threshold | $\varepsilon_{\Theta}$ | $10^{-6}$ s |
| $z$-axis Aerodynamic Coefficient | $C_z$ | $7 \pm 3.5$ | Integral Parameter 2 Threshold | $\varepsilon_{\tilde{\Theta}}$ | $10^{-6}$ s |
| Thrust Ratio | $n$ | 3 | Terminal Height Threshold | $\varepsilon_{yt}$ | 0 m |
| Time Allocation Ratio | $\mu$ | 0.5 | Time Allocation Ratio Threshold | $\varepsilon_{\mu}$ | $10^{-6}$ |
| Guidance Cycle | $\Delta\tau$ | 10 ms | Remaining Time Threshold | $\varepsilon_t$ | 0.1 s |

The Monte Carlo simulation results are shown in Fig. 11.

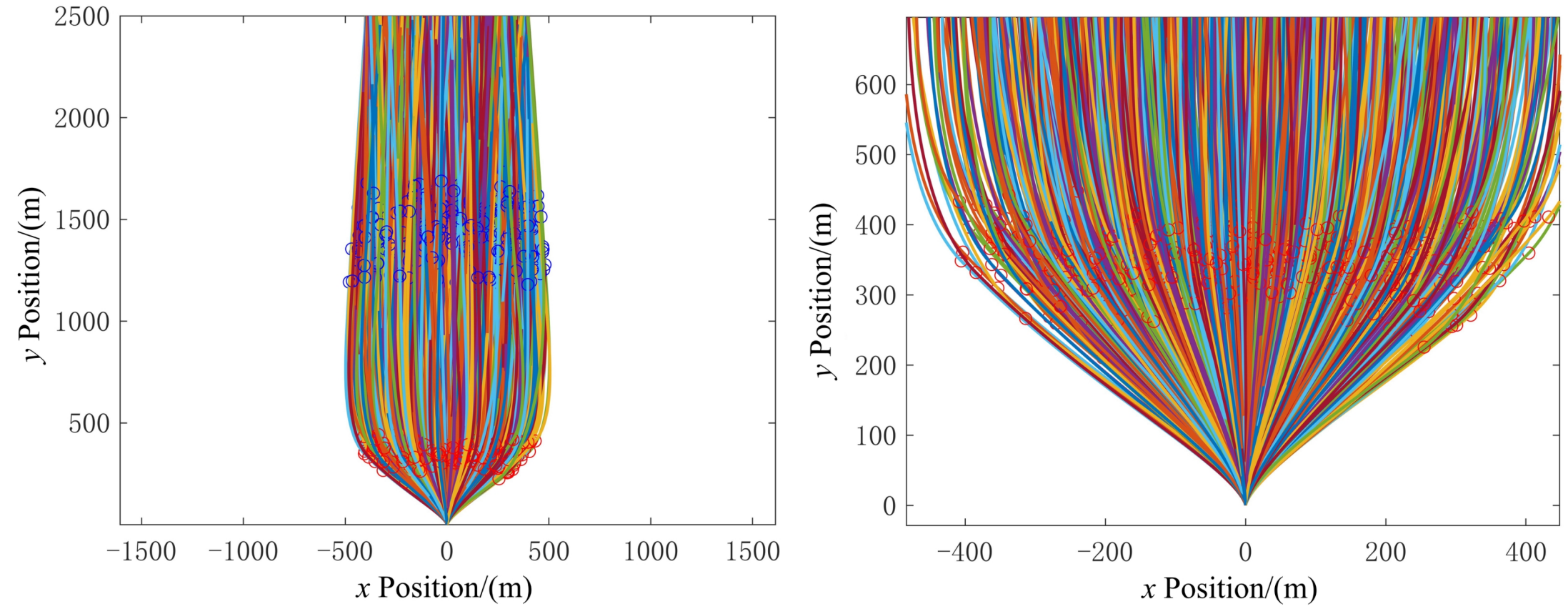

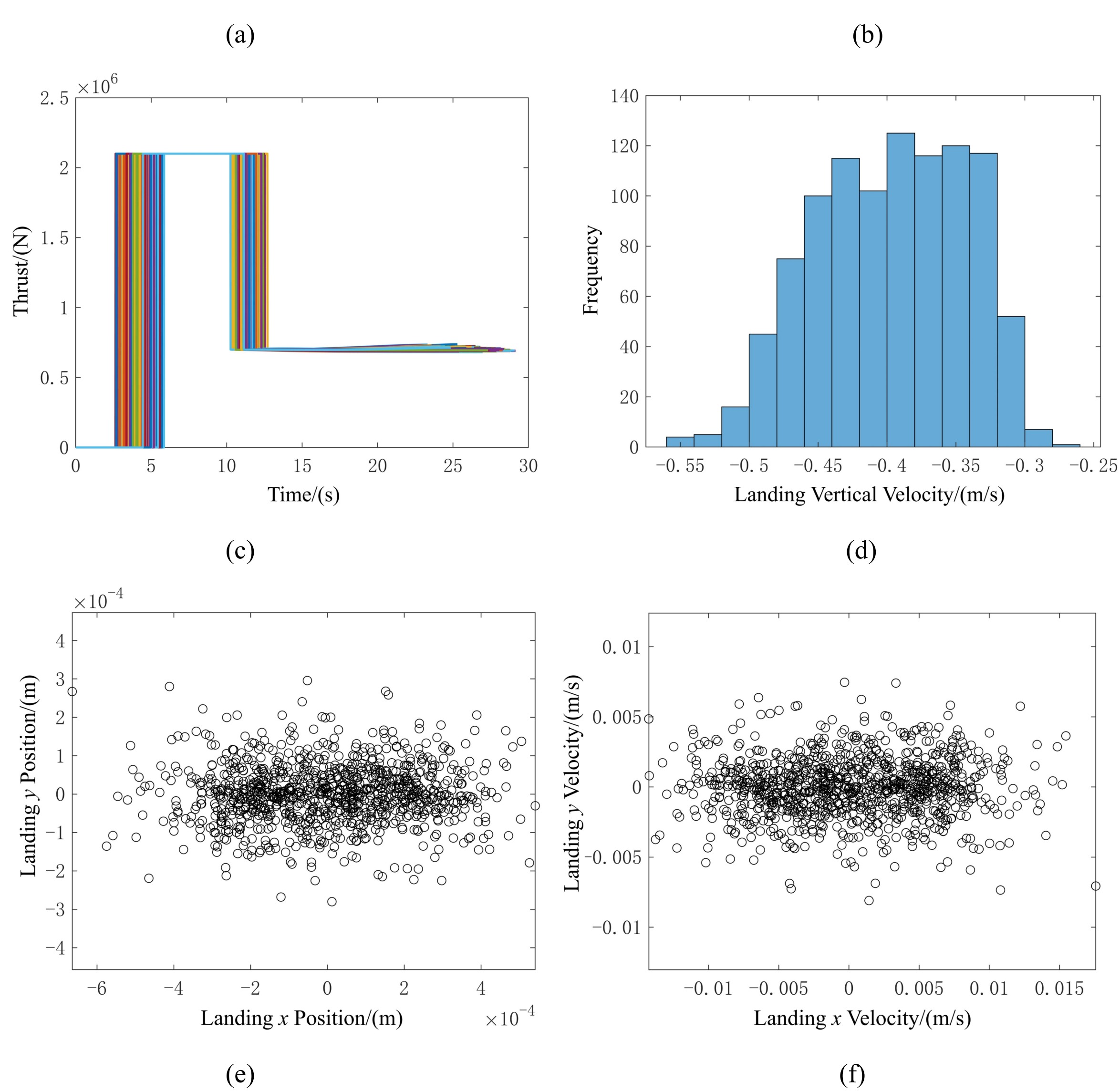

**Fig. 11 Monte Carlo simulation results. (a) Trajectory curves. (b) Locally magnified trajectory curves. (c) Thrust magnitude vs. time. (d) Landing vertical velocity distribution. (e) Landing horizontal position distribution. (f) Landing horizontal velocity distribution.**

Fig. 11 (a) shows the typical flight trajectories of the rocket during the landing. The blue circles mark the positions of engine activation, and the red circles mark the positions of thrust switching. The results indicate that the proposed guidance algorithm can adaptively determine the timing of guidance activation and thrust switching based on the real-time flight state, and all simulated trajectories stably converge to the target landing point (the origin of coordinates).

As seen from the thrust curves in Fig. 11 (c), the engine activation times are distributed between 2.6 s and 5.6 s. After activation, the engine first enters the high-thrust phase, maintaining precisely at 3 times the nominal thrust (2100 kN) to achieve rapid deceleration; it then switches to the low-thrust phase between 10.3 s and 12.6 s, with the thrust dropping to 1 times the nominal value. Affected by system parameter deviations and disturbances, the thrust command exhibits slight fluctuations, with a maximum value of 737.608 kN and a minimum value of 681.293 kN. The engine only needs to have a thrust adjustment range of 97% to 105% to avoid thrust saturation. This range is only 8.04% of the nominal thrust, significantly reducing the requirements for the adjustment capability of the propulsion system. This is mainly attributed to the adopted dual-thrust mode: the rocket achieves rapid deceleration by flexibly adjusting the duration of the high-thrust phase; after entering the low-thrust phase, the aerodynamic influence is greatly weakened, and at this point, only fine-tuning of the thrust is required to achieve landing.

Fig. 11 (d) shows that the vertical velocity at landing is concentrated between 0.25 m/s and 0.55 m/s, satisfying the buffering constraints for safe landing. Furthermore, Fig. 11 (e) and (f) show the landing accuracy in the horizontal direction: the horizontal position error is less than 0.6 mm, and the horizontal velocity error is below 1.5 cm/s. Both are at extremely low levels, indicating that the proposed method possesses excellent terminal accuracy. This performance also benefits from the application of the dual-thrust mode: in the low-thrust phase, the thrust-to-weight ratio can be reduced to a relatively low level (less than 1.5), thereby creating conditions for landing with centimeter-level accuracy.

## VIII. Conclusions

This paper proposes a dual-thrust mode analytical guidance method for powered landing (SAM-DS) that aims to achieve optimal attitude angle smoothness. Through theoretical derivation, it is proven that the optimal attitude command is a piecewise cubic function, thereby transforming the trajectory optimization problem into a parametric analytical optimization problem. The proposed method adopts a three-phase framework ("Pre-activation – High-thrust phase – Low-thrust phase") and incorporates an aerodynamic correction algorithm to achieve high-precision landing in complex environments.

The research results indicate that the algorithm significantly reduces the demand for engine thrust modulation capability (8%) while ensuring extremely high computational efficiency (less than 0.6 ms) and optimality. Under conditions of severe parametric deviations, the algorithm still achieves high-precision fixed-point landing with terminal accuracy satisfying engineering safety standards, demonstrating excellent robustness and engineering applicability. The proposed framework features clear logic and strong interpretability, providing an efficient and reliable solution for reusable rocket guidance systems.

## Appendix A: Proof of Optimality

**Proof:** First, the Hamiltonian function is constructed. If $\frac{(n-1)T}{I_{\mathrm{sp}}m_0}\mu t_{\mathrm{go}}^{(2)}+\frac{T}{I_{\mathrm{sp}}m_0}\tau \ll 1$, $\left|\frac{1}{2}\cos\bar{\theta}\Delta\theta^2\right| \ll 1$ and $\left|\frac{1}{2}\sin\bar{\theta}\Delta\theta^2\right| \ll 1$, the Hamiltonian function can be simplified as:

$$H=\begin{cases}\frac{1}{2}\left(u_1^2+u_2^2\right)+\lambda_{rx}v_x+\lambda_{ry}v_y+\lambda_{rz}v_z \\ +\lambda_{vx}\left[\frac{nT}{m_0}\left(\cos\bar{\theta}-\sin\bar{\theta}\Delta\theta\right)\right]+\lambda_{vy}\left[\frac{nT}{m_0}\left(\sin\bar{\theta}+\cos\bar{\theta}\Delta\theta\right)\right] \\ +\lambda_{vz}\left[-\frac{nT}{m_0}\left(\bar{\psi}+\Delta\psi\right)\right]+\lambda_\theta u_1+\lambda_\psi u_2+\lambda_T u_3+\ldots\ ,\ \tau\le\mu t_{\mathrm{go}}^{(2)} \\ \frac{1}{2}\left(u_1^2+u_2^2\right)+\lambda_{rx}v_x+\lambda_{ry}v_y+\lambda_{rz}v_z \\ +\lambda_{vx}\left[\frac{T}{m_0}\left(\cos\bar{\theta}-\sin\bar{\theta}\Delta\theta\right)\right]+\lambda_{vy}\left[\frac{T}{m_0}\left(\sin\bar{\theta}+\cos\bar{\theta}\Delta\theta\right)\right] \\ +\lambda_{vz}\left[-\frac{T}{m_0}\left(\bar{\psi}+\Delta\psi\right)\right]+\lambda_\theta u_1+\lambda_\psi u_2+\lambda_T u_3+\ldots\ ,\ \mu t_{\mathrm{go}}^{(2)}<\tau\le(\mu+1)t_{\mathrm{go}}^{(2)}\end{cases} \tag{A1}$$

According to the Minimum Principle, it follows that:

$$\frac{\partial H}{\partial u_1}=u_1+\lambda_\theta,\frac{\partial H}{\partial u_2}=u_2+\lambda_\psi,\frac{\partial H}{\partial u_3}=\lambda_T \tag{A2}$$

Considering the constraint $u_3=0$, the optimal solution to the problem can be expressed in terms of the costate variables as:

$$u_1^*=-\lambda_\theta,u_2^*=-\lambda_\psi,u_3^*=0 \tag{A3}$$

The costate equations are listed as follows:

$$\begin{aligned}&\dot{\lambda}_\theta=\begin{cases}\lambda_{vx}\frac{nT}{m_0}\sin\bar{\theta}-\lambda_{vy}\frac{nT}{m_0}\cos\bar{\theta},\ \tau\le\mu t_{\mathrm{go}}^{(2)} \\ \lambda_{vx}\frac{T}{m_0}\sin\bar{\theta}-\lambda_{vy}\frac{T}{m_0}\cos\bar{\theta},\ \mu t_{\mathrm{go}}^{(2)}<\tau\le(\mu+1)t_{\mathrm{go}}^{(2)}\end{cases} \\ &\dot{\lambda}_\psi=\begin{cases}\lambda_{vz}\frac{nT}{m_0},\ \tau\le\mu t_{\mathrm{go}}^{(2)} \\ \lambda_{vz}\frac{T}{m_0},\ \mu t_{\mathrm{go}}^{(2)}<\tau\le(\mu+1)t_{\mathrm{go}}^{(2)}\end{cases} \\ &\dot{\lambda}_{vx}=-\lambda_{rx},\dot{\lambda}_{vy}=-\lambda_{ry},\dot{\lambda}_{vz}=-\lambda_{rz} \\ &\dot{\lambda}_{rx}=0,\dot{\lambda}_{ry}=0,\dot{\lambda}_{rz}=0\end{aligned} \tag{A4}$$

From Eq. (A4), it is known that $\lambda_{rx},\lambda_{ry},\lambda_{rz}$ are constant functions. Without loss of generality, let:

$$\lambda_{rx}=\lambda_{rx0},\lambda_{ry}=\lambda_{ry0},\lambda_{rz}=\lambda_{rz0} \tag{A5}$$

From $\dot{\lambda}_{vx}=-\lambda_{rx},\dot{\lambda}_{vy}=-\lambda_{ry},\dot{\lambda}_{vz}=-\lambda_{rz}$, it is known that $\lambda_{vx},\lambda_{vy},\lambda_{vz}$ are linear functions:

$$\lambda_{vx}=-\lambda_{rx0}\tau+\lambda_{vx0},\lambda_{vy}=-\lambda_{ry0}\tau+\lambda_{vy0},\lambda_{vz}=-\lambda_{rz0}\tau+\lambda_{vz0} \tag{A6}$$

Then, $\dot{\lambda}_\theta$ and $\dot{\lambda}_\psi$ can be rearranged into the following forms according to (A4):

$$\begin{aligned}
\dot{\lambda}_\theta &= \begin{cases} \dfrac{nT}{m_0}\left(-\lambda_{rx0}\sin\bar{\theta}+\lambda_{ry0}\cos\bar{\theta}\right)\tau+\dfrac{nT}{m_0}\left(\lambda_{vx0}\sin\bar{\theta}-\lambda_{vy0}\cos\bar{\theta}\right),\ \tau\le\mu t_{\mathrm{go}}^{(2)} \\ \dfrac{T}{m_0}\left(-\lambda_{rx0}\sin\bar{\theta}+\lambda_{ry0}\cos\bar{\theta}\right)\tau+\dfrac{T}{m_0}\left(\lambda_{vx0}\sin\bar{\theta}-\lambda_{vy0}\cos\bar{\theta}\right),\ \mu t_{\mathrm{go}}^{(2)}<\tau\le(\mu+1)t_{\mathrm{go}}^{(2)} \end{cases} \\
\dot{\lambda}_\psi &= \begin{cases} \dfrac{nT}{m_0}\left(-\lambda_{rz0}\tau+\lambda_{vz0}\right),\ \tau\le\mu t_{\mathrm{go}}^{(2)} \\ \dfrac{T}{m_0}\left(-\lambda_{rz0}\tau+\lambda_{vz0}\right),\ \mu t_{\mathrm{go}}^{(2)}<\tau\le(\mu+1)t_{\mathrm{go}}^{(2)} \end{cases}
\end{aligned} \tag{A7}$$

Obviously, both $\dot{\lambda}_\theta$ and $\dot{\lambda}_\psi$ are piecewise linear functions. Further integration reveals that $\lambda_\theta$ and $\lambda_\psi$ are piecewise quadratic functions, satisfying:

$$\lambda_\theta\Big|_{\tau=t_{\mathrm{go}}^{(1)-}}=\lambda_\theta\Big|_{\tau=t_{\mathrm{go}}^{(1)+}},\lambda_\psi\Big|_{\tau=t_{\mathrm{go}}^{(1)-}}=\lambda_\psi\Big|_{\tau=t_{\mathrm{go}}^{(1)+}} \tag{A8}$$

Further, based on Eq. (A3) and the integral of $u_1,u_2$, it can be concluded that $\Delta\theta^*$ and $\Delta\psi^*$ are piecewise cubic functions, satisfying:

$$\Delta\theta^*\Big|_{\tau=t_{\mathrm{go}}^{(1)-}}=\Delta\theta^*\Big|_{\tau=t_{\mathrm{go}}^{(1)+}},\Delta\psi^*\Big|_{\tau=t_{\mathrm{go}}^{(1)-}}=\Delta\psi^*\Big|_{\tau=t_{\mathrm{go}}^{(1)+}} \tag{A9}$$

At this point, it has been proven that the optimal solutions $\Delta\theta^*$ and $\Delta\psi^*$ to Problem 1 are cubic functions of time. The continuity conditions can be derived from Eq. (A9), the differentiability conditions from Eq. (A8), and the proportionality conditions from Eq. (A7).

## Appendix B: Representation of Integral Functionals

First, variables and functions are defined:

$$\begin{aligned}
f^{(1)}(x) &= \frac{nt_{\mathrm{go}}^{(2)}}{(1+x)(2+x)}\left\{\begin{aligned}&(1+\mu)\left[-t_{\mathrm{go}}^{(2)}(1+\mu)\right]^x\left[2+x+hnt_{\mathrm{go}}^{(2)}(1+\mu)\right] \\ &-\left(-t_{\mathrm{go}}^{(2)}\right)^x\left\{2+x+hnt_{\mathrm{go}}^{(2)}\left[1+(2+x)\mu\right]\right\}\end{aligned}\right\} \\
f^{(2)}(x) &= \frac{t_{\mathrm{go}}^{(2)}}{(1+x)(2+x)}\left(-t_{\mathrm{go}}^{(2)}\right)^x\left\{2+x+ht_{\mathrm{go}}^{(2)}\left[1+n(2+x)\mu\right]\right\}
\end{aligned} \tag{B1}$$

To simplify the formula expressions, the following vectors are defined:

$$\begin{aligned}
\boldsymbol{p}^{(i)} &= \left[d^{(i)},c^{(i)},b^{(i)},a^{(i)}\right]^T, i=1,2 \\
\boldsymbol{p}'^{(i)} &= \left[d'^{(i)},c'^{(i)},b'^{(i)},a'^{(i)}\right]^T, i=1,2
\end{aligned} \tag{B2}$$

Assume the pitch angle profile is given by Eq. (25) and the yaw angle profile by Eq. (37). Then, the integral functional defined in Eq. (10) can be expanded and simplified sequentially. First, calculate $\Theta(\Delta\theta)$:

$$
\begin{aligned}
\Theta(\Delta\theta) &= n\Theta^{(1)}(\Delta\theta)+\Theta^{(2)}(\Delta\theta) \\
&= n\int_{\tau=0}^{\mu t_{\mathrm{go}}^{(2)}} \Delta\theta(1+nh\tau)\,\mathrm{d}\tau + \int_{\tau=\mu t_{\mathrm{go}}^{(2)}}^{(\mu+1)t_{\mathrm{go}}^{(2)}} \Delta\theta\left[1+(n-1)h\mu t_{\mathrm{go}}^{(2)}+h\tau\right]\mathrm{d}\tau \\
&= \sum_{i=1}^{2}\sum_{j=1}^{4} p_j^{(i)} f^{(i)}(j-1)
\end{aligned} \tag{B3}
$$

Substituting Eq. (25) into Eq. (B3) to eliminate $a^{(1)} \sim d^{(1)}$ yields:

$$
\begin{aligned}
\Theta(\Delta\theta) &= na^{(2)} f^{(1)}(3) + nb^{(2)} f^{(1)}(2) + \left[c^{(2)} + (n-1)t_{\mathrm{go}}^{(2)}\left(2b^{(2)} - 3a^{(2)}t_{\mathrm{go}}^{(2)}\right)\right] f^{(1)}(1) \\
&+\left[d^{(2)} + (n-1)t_{\mathrm{go}}^{(2)2}\left(b^{(2)} - 2a^{(2)}t_{\mathrm{go}}^{(2)}\right)\right] f^{(1)}(0) \\
&+a^{(2)} f^{(2)}(3) + b^{(2)} f^{(2)}(2) + c^{(2)} f^{(2)}(1) + d^{(2)} f^{(2)}(0)
\end{aligned} \tag{B4}
$$

After rearranging terms and simplifying, we obtain:

$$
\Theta(\Delta\theta) = A_1 a^{(2)} + B_1 b^{(2)} + C_1 c^{(2)} + D_1 d^{(2)} \tag{B5}
$$

where

$$
\begin{aligned}
A_1 &= nf^{(1)}(3) + f^{(2)}(3) - 3(n-1)t_{\mathrm{go}}\ f^{(1)}(1) - 2(n-1)t_{\mathrm{go}}^{2}\ f^{(1)}(0) \\
B_1 &= nf^{(1)}(2) + f^{(2)}(2) + 2(n-1)t_{\mathrm{go}}^{(2)} f^{(1)}(1) + (n-1)t_{\mathrm{go}}^{(2)2} f^{(1)}(0) \\
C_1 &= f^{(1)}(1) + f^{(2)}(1) \\
D_1 &= f^{(1)}(0) + f^{(2)}(0)
\end{aligned} \tag{B6}
$$

The yaw channel is analogous to the pitch channel, and the result can be written directly as:

$$
\begin{aligned}
\Theta(\Delta\psi) &= n\Theta^{(1)}(\Delta\psi)+\Theta^{(2)}(\Delta\psi) \\
&= n\int_{\tau=0}^{\mu t_{\mathrm{go}}^{(2)}} \Delta\psi(1+nh\tau)\,\mathrm{d}\tau + \int_{\tau=\mu t_{\mathrm{go}}^{(2)}}^{(\mu+1)t_{\mathrm{go}}^{(2)}} \Delta\psi\left[1+(n-1)h\mu t_{\mathrm{go}}^{(2)}+h\tau\right]\mathrm{d}\tau \\
&= A_1 a'^{(2)} + B_1 b'^{(2)} + C_1 c'^{(2)} + D_1 d'^{(2)}
\end{aligned} \tag{B7}
$$

In the aerodynamic correction section, it is necessary to calculate the velocity magnitude at the switching point, which requires $n\Theta^{(1)}(\Delta\theta)$ and $n\Theta^{(1)}(\Delta\psi)$. These can be directly written according to Eqs. (B3) and (B7) as follows:

$$
\begin{aligned}
n\Theta^{(1)}(\Delta\theta) &= n\int_{\tau=0}^{\mu t_{\mathrm{go}}^{(2)}} \Delta\theta(1+nh\tau)\,\mathrm{d}\tau = \sum_{j=1}^{4} p_j^{(1)} f^{(1)}(j-1) \\
n\Theta^{(1)}(\Delta\psi) &= n\int_{\tau=0}^{\mu t_{\mathrm{go}}^{(2)}} \Delta\psi(1+nh\tau)\,\mathrm{d}\tau = \sum_{j=1}^{4} p_j'^{(1)} f^{(1)}(j-1)
\end{aligned} \tag{B8}
$$

Following the same approach as the derivation of $\Theta(\Delta\theta)$, calculate the double integral functional $\tilde{\Theta}(\Delta\theta)$:

$$
\begin{aligned}
\tilde{\Theta}(\Delta\theta) &= nt_{\mathrm{go}}^{(2)}\Theta^{(1)}(\Delta\theta) + n\tilde{\Theta}^{(1)}(\Delta\theta) + \tilde{\Theta}^{(2)}(\Delta\theta) \\
&= nt_{\mathrm{go}}^{(2)} \int_{\tau=0}^{\mu t_{\mathrm{go}}^{(2)}} \Delta\theta(1+nh\tau)\,\mathrm{d}\tau + n\int_{t=0}^{\mu t_{\mathrm{go}}^{(2)}}\int_{\tau=0}^{t} \Delta\theta(1+nh\tau)\,\mathrm{d}\tau\mathrm{d}t \\
&+ \int_{t=\mu t_{\mathrm{go}}^{(2)}}^{(\mu+1)t_{\mathrm{go}}^{(2)}}\int_{\tau=\mu t_{\mathrm{go}}^{(2)}}^{t} \Delta\theta\left[1+(n-1)h\mu t_{\mathrm{go}}^{(2)} + h\tau\right]\mathrm{d}\tau\mathrm{d}t \\
&= -\sum_{i=1}^{2}\sum_{j=1}^{4} p_j^{(i)} f^{(i)}(j)
\end{aligned} \tag{B9}
$$

Substituting Eq. (25) into Eq. (B9) to eliminate $a^{(1)} \sim d^{(1)}$ yields:

$$
\begin{aligned}
\tilde{\Theta}(\Delta\theta) &= -na^{(2)}f^{(1)}(4) - nb^{(2)}f^{(1)}(3) - \left[c^{(2)} + (n-1)t_{\mathrm{go}}^{(2)}\left(2b^{(2)} - 3a^{(2)}t_{\mathrm{go}}^{(2)}\right)\right]f^{(1)}(2) \\
&-\left[d^{(2)} + (n-1)t_{\mathrm{go}}^{(2)2}\left(b^{(2)} - 2a^{(2)}t_{\mathrm{go}}^{(2)}\right)\right]f^{(1)}(1) \\
&-a^{(2)}f^{(2)}(4) - b^{(2)}f^{(2)}(3) - c^{(2)}f^{(2)}(2) - d^{(2)}f^{(2)}(1)
\end{aligned} \tag{B10}
$$

After rearranging terms and simplifying, we obtain:

$$
\tilde{\Theta}(\Delta\theta) = A_2a^{(2)} + B_2b^{(2)} + C_2c^{(2)} + D_2d^{(2)} \tag{B11}
$$

where

$$
\begin{aligned}
A_2 &= -\left[nf^{(1)}(4) + f^{(2)}(4) - 3(n-1)t_{\mathrm{go}}^{(2)2}f^{(1)}(2) - 2(n-1)t_{\mathrm{go}}^{(2)3}f^{(1)}(1)\right] \\
B_2 &= -\left[nf^{(1)}(3) + f^{(2)}(3) + 2(n-1)t_{\mathrm{go}}^{(2)}f^{(1)}(2) + (n-1)t_{\mathrm{go}}^{(2)2}f^{(1)}(1)\right] \\
C_2 &= -\left[f^{(1)}(2) + f^{(2)}(2)\right] \\
D_2 &= -\left[f^{(1)}(1) + f^{(2)}(1)\right]
\end{aligned} \tag{B12}
$$

The yaw channel is analogous to the pitch channel, and the result can be written directly as:

$$
\begin{aligned}
\tilde{\Theta}(\Delta\psi) &= nt_{\mathrm{go}}^{(2)}\Theta^{(1)}(\Delta\psi) + n\tilde{\Theta}^{(1)}(\Delta\psi) + \tilde{\Theta}^{(2)}(\Delta\psi) \\
&= nt_{\mathrm{go}}^{(2)} \int_{\tau=0}^{\mu t_{\mathrm{go}}^{(2)}} \Delta\psi(1+nh\tau)\,\mathrm{d}\tau + n\int_{t=0}^{\mu t_{\mathrm{go}}^{(2)}}\int_{\tau=0}^{t} \Delta\psi(1+nh\tau)\,\mathrm{d}\tau\mathrm{d}t \\
&+ \int_{t=\mu t_{\mathrm{go}}^{(2)}}^{(\mu+1)t_{\mathrm{go}}^{(2)}}\int_{\tau=\mu t_{\mathrm{go}}^{(2)}}^{t} \Delta\psi\left[1+(n-1)h\mu t_{\mathrm{go}}^{(2)} + h\tau\right]\mathrm{d}\tau\mathrm{d}t \\
&= A_2a'^{(2)} + B_2b'^{(2)} + C_2c'^{(2)} + D_2d'^{(2)}
\end{aligned} \tag{B13}
$$

Finally, derive the integral functionals $\Theta(\Delta\theta^2)$ and $\tilde{\Theta}(\Delta\theta^2)$:

$$
\begin{aligned}
\Theta(\Delta\theta^2) &= n\Theta^{(1)}(\Delta\theta^2) + \Theta^{(2)}(\Delta\theta^2) \\
&= n\int_{\tau=0}^{\mu t_{\mathrm{go}}^{(2)}} \Delta\theta^2(1+nh\tau)\,\mathrm{d}\tau \\
&+ \int_{\tau=\mu t_{\mathrm{go}}^{(2)}}^{(\mu+1)t_{\mathrm{go}}^{(2)}} \Delta\theta^2\left(1+(n-1)h\mu t_{\mathrm{go}}^{(2)} + h\tau\right)\mathrm{d}\tau \\
&= \sum_{i=1}^{2}\sum_{j=1}^{4}\sum_{k=1}^{4} p_j^{(i)}p_k^{(i)}f^{(i)}(j+k-2) \\
&= \sum_{i=1}^{2}\left[\sum_{j=1}^{4} p_j^{(i)2}f^{(i)}(2j-2) + 2\sum_{1\le j<k\le 4} p_j^{(i)}p_k^{(i)}f^{(i)}(j+k-2)\right]
\end{aligned} \tag{B14}
$$

$$
\begin{aligned}
\tilde{\Theta}\left(\Delta\theta^2\right) &= n t_{\text{go}}^{(2)} \Theta^{(1)}\left(\Delta\theta^2\right) + n\tilde{\Theta}^{(1)}\left(\Delta\theta^2\right) + \tilde{\Theta}^{(2)}\left(\Delta\theta^2\right) \\
&= n t_{\text{go}}^{(2)} \int_{\tau=0}^{\mu t_{\text{go}}^{(2)}} \Delta\theta^2 \left(1+nh\tau\right) \mathrm{d}\tau + n \int_{t=0}^{\mu t_{\text{go}}^{(2)}} \int_{\tau=0}^{t} \Delta\theta^2 \left(1+nh\tau\right) \mathrm{d}\tau \mathrm{d}t \\
&+ \int_{t=\mu t_{\text{go}}^{(2)}}^{(\mu+1) t_{\text{go}}^{(2)}} \int_{\tau=\mu t_{\text{go}}^{(2)}}^{t} \Delta\theta^2 \left[1+\left(n-1\right) h \mu t_{\text{go}}^{(2)} + h\tau\right] \mathrm{d}\tau \mathrm{d}t \\
&= -\sum_{i=1}^{2} \sum_{j=1}^{4} \sum_{k=1}^{4} p_j^{(i)} p_k^{(i)} f^{(i)} \left(j+k-1\right) \\
&= -\sum_{i=1}^{2} \left[ \sum_{j=1}^{4} p_j^{(i)2} f^{(i)} \left(2j-1\right) + 2 \sum_{1\le j<k\le 4} p_j^{(i)} p_k^{(i)} f^{(i)} \left(j+k-1\right) \right]
\end{aligned} \tag{B15}
$$

In the aerodynamic correction section, it is necessary to calculate the velocity at the switching point, which requires $n\Theta^{(1)}\left(\Delta\theta^2\right)$. This can be written according to Eq. (B14) as follows:

$$
\begin{aligned}
n\Theta^{(1)}\left(\Delta\theta^2\right) &= n \int_{\tau=0}^{\mu t_{\text{go}}^{(2)}} \Delta\theta^2 \left(1+nh\tau\right) \mathrm{d}\tau \\
&= \sum_{j=1}^{4} p_j^{(1)2} f^{(1)} \left(2j-2\right) + 2 \sum_{1\le j<k\le 4} p_j^{(1)} p_k^{(1)} f^{(1)} \left(j+k-2\right)
\end{aligned} \tag{B16}
$$

## Appendix C: Solution to the Remaining Time Equation

From Eq. (22), the specific form of the remaining flight time equation is as follows:

$$
\begin{aligned}
&\tilde{A} t_{\text{go}}^{(2)4} + \tilde{B} t_{\text{go}}^{(2)3} + \tilde{C} t_{\text{go}}^{(2)2} + \tilde{D} t_{\text{go}}^{(2)} + \tilde{E} = 0 \\
&\tilde{A} = \frac{\left(n\mu+1\right)^4 T^4}{4 I_{\text{sp}}^2 m_0^4} \\
&\tilde{B} = \frac{\left(n\mu+1\right)^3 T^3}{I_{\text{sp}} m_0^3} \\
&\tilde{C} = \left(n\mu+1\right)^2 \frac{T^2}{m_0^2} - \frac{\left(n\mu+1\right)^2 T^3}{2 I_{\text{sp}} m_0^3} \Theta\left(\Delta\theta^2\right) - \left[g\left(\mu+1\right) + k \left\|\boldsymbol{v}_0\right\| \Pi v_{y0}\right]^2 - k^2 \left\|\boldsymbol{v}_0\right\|^2 \Pi^2 v_{x0}^2 \\
&\tilde{D} = -\frac{\left(n\mu+1\right) T^2}{m_0^2} \Theta\left(\Delta\theta^2\right) + 2\left[g\left(\mu+1\right) + k \left\|\boldsymbol{v}_0\right\| \Pi v_{y0}\right]\left(v_{y0}+\xi_y\right) + 2k \left\|\boldsymbol{v}_0\right\| \Pi v_{x0} \left(v_{x0}+\xi_x\right) \\
&\tilde{E} = \frac{T^2}{4 m_0^2} \Theta\left(\Delta\theta^2\right)^2 - \left(v_{y0}+\xi_y\right)^2 - \left(v_{x0}+\xi_x\right)^2
\end{aligned} \tag{C1}
$$

The above quartic equation can be solved according to the following process. A quartic equation has four solutions, and the unique physically meaningful solution is:

$$
\begin{aligned}
&\Lambda = 2\tilde{C}^3 - 9\tilde{B}\tilde{C}\tilde{D} + 27\tilde{A}\tilde{D}^2 + 27\tilde{B}^2\tilde{E} - 72\tilde{A}\tilde{C}\tilde{E} \\
&\Sigma = \tilde{C}^2 - 3\tilde{B}\tilde{D} + 12\tilde{A}\tilde{E} \\
&\Xi = 4\Sigma^3 - \Lambda^2 \\
&\Gamma = \frac{1}{3}\arctan\left(\frac{\sqrt{\Xi}}{\Lambda}\right) \\
&\Omega = \frac{\tilde{B}^2}{4\tilde{A}^2} - \frac{2\tilde{C}}{3\tilde{A}} + \frac{2}{3\tilde{A}}\Sigma^{\frac{1}{2}}\cos\Gamma \\
&t_{\mathrm{go}}^{(2)} = -\frac{\tilde{B}}{4\tilde{A}} + \frac{1}{2}\Omega^{\frac{1}{2}} + \frac{1}{2}\left[\frac{3\tilde{B}^2}{4\tilde{A}^2} - 2\frac{\tilde{C}}{\tilde{A}} - \Omega + \left(-\frac{\tilde{B}^3}{4\tilde{A}^3} + \frac{\tilde{B}\tilde{C}}{\tilde{A}^2} - \frac{2\tilde{D}}{\tilde{A}}\right)\Omega^{-\frac{1}{2}}\right]^{\frac{1}{2}}
\end{aligned}
\tag{C2}
$$

The above solution method does not require iteration, offers superior computational efficiency, and eliminates the need to consider convergence issues.